\documentclass[12pt,twoside]{amsart}

\usepackage{amsmath}
\usepackage{amssymb}
\usepackage{amscd}

\DeclareMathSymbol{\ordcol}{\mathord}{operators}{'072}
\newcommand{\col}{{{\hskip1.5pt\ordcol\hskip1.5pt}}}

\renewcommand{\bar}{\overline}

\newcommand{\eps}{\varepsilon}

\newcommand{\CC}{\mathbb{C}}
\newcommand{\DD}{\mathbb{D}}

\newcommand{\FF}{\mathbb{F}}

\newcommand{\PP}{\mathbb{P}}
\newcommand{\QQ}{\mathbb{Q}}
\newcommand{\RR}{\mathbb{R}}

\newcommand{\ZZ}{\mathbb{Z}}

\newcommand{\Cv}{\CC_v}

\newcommand{\ints}{{\mathcal O}}

\newcommand{\maxid}{{\mathcal M}}

\newcommand{\calK}{{\mathcal K}}

\newcommand{\calP}{{\mathcal P}}

\newcommand{\PCv}{\PP^1(\Cv)}

\newcommand{\PK}{\PP^1(K)}

\newcommand{\PGL}{\hbox{\rm PGL}}

\newcommand{\dist}{\hbox{\rm dist}}

\newcommand{\Dbar}{\bar{D}}

\newcommand{\dsps}{\displaystyle}

\theoremstyle{plain}
\newtheorem{thm}{Theorem}[section]

\newtheorem{lemma}[thm]{Lemma}
\newtheorem{defin}[thm]{Definition}

\newtheorem*{mthm}{Main Theorem}

\newtheorem*{ufbconj}{Uniform Boundedness Conjecture}

\theoremstyle{definition}

\newtheorem{remark}[thm]{Remark}

\newtheorem{example}[thm]{Example}

\title[Preperiodic points]{Preperiodic points
  of polynomials over global fields}
\author{Robert L. Benedetto}
\date{June 21, 2005}
\thanks{The author gratefully acknowledges the support
  of a Miner D.\ Crary Research Fellowship from Amherst College
  and an NSA Young Investigator Grant}
\address{Department of Mathematics and Computer Science \\
        Amherst College \\
        Amherst, MA 01002 \\
        USA}
\email{rlb@cs.amherst.edu}
\urladdr{http://www.cs.amherst.edu/\textasciitilde rlb}

\begin{document}

\newcounter{bean}
\newcounter{sheep}

\begin{abstract}
Given a global field $K$ and a polynomial $\phi$ defined
over $K$ of degree at least two, Morton and Silverman
conjectured in 1994 that the number of $K$-rational
preperiodic points of $\phi$ is bounded
in terms of only the degree of $K$ and the degree of $\phi$.
In 1997, for quadratic polynomials over $K=\QQ$,
Call and Goldstine proved a bound which was exponential in $s$,
the number of primes of bad reduction of $\phi$.
By careful analysis
of the filled Julia sets at each prime,
we present an improved bound on the order of $s\log s$.
Our bound applies to polynomials of any degree (at least two)
over any global field $K$.
\end{abstract}

\maketitle

Let $K$ be a field, and let $\phi\in K(z)$ be
a rational function.  Let $\phi^n$ denote the $n^{\text{th}}$
iterate of $\phi$ under composition; that is,
$\phi^0$ is the identity function, and for $n\geq 1$,
$\phi^n = \phi\circ \phi^{n-1}$.
We will study the dynamics $\phi$ on the projective line
$\PP^1(K)$.  In particular, we say a point $x$ is
{\em preperiodic} under $\phi$ if there are integers $n>m\geq 0$
such that $\phi^m(x)=\phi^n(x)$.  The point
$y=\phi^m(x)$ satisfies $\phi^{n-m}(y)=y$
and is said to be {\em periodic}, as its iterates
will forever cycle through the same finite sequence of values.
Note that $x\in\PP^1(K)$ is preperiodic if and only
if its orbit $\{\phi^n(x):n\geq 0\}$ is finite.

For example, let $K=\QQ$ and $\phi(z)=z^2 - 29/16$.  Then
$\{5/4, -1/4, -7/4\}$ forms a periodic cycle
(of period $3$), and $-5/4$, $1/4$, $7/4$, and $\pm 3/4$ each
land on this cycle after one or two iterations.  In addition,
the point $\infty$ is of course fixed.  These nine $\QQ$-rational
points are all preperiodic under $\phi$.
Meanwhile, it is not difficult to see that no other
point in $\PP^1(\QQ)$ is preperiodic by showing that
the denominator of a rational preperiodic point must be $4$,
and that the absolute value must be less than $2$.

In general, for any global field $K$,
any dimension $N\geq 1$, and any morphism
$\phi:\PP^N\rightarrow\PP^N$ over $K$ of degree
at least two, Northcott proved in 1950 that the number of
$K$-rational preperiodic points of $\phi$ is finite \cite{Nor}.
More precisely, he showed that the
preperiodic points form a set of bounded arithmetic height.
Years later,
by analogy with the Theorems of Mazur \cite{Maz} and Merel \cite{Mer}
on $K$-rational torsion of elliptic curves,
Morton and Silverman proposed the following Conjecture \cite{MS1}.
\begin{ufbconj}
{\rm (Morton and Silverman, 1994)} \\
Given integers $D,N\geq 1$ and $d\geq 2$, there is a constant
$\kappa=\kappa(D,N,d)$ with the following property.
Let $K$ be a number field with $[K\col\QQ \hspace{0.75pt}]=D$,
and let $\phi:\PP^N\rightarrow\PP^N$ be a morphism of degree $d$
defined over $K$.
Then $\phi$ has at most $\kappa$ preperiodic
points in $\PP^N(K)$.
\end{ufbconj}
The analogy between preperiodic points and torsion
comes from the fact that
the torsion points of an elliptic curve $E$
are precisely the preperiodic points of the
multiplication-by-two map $[2]:E\rightarrow E$.
In fact, taking $x$-coordinates, the map $[2]$
induces a rational function (known to dynamicists
as a Latt\`{e}s map)
$\phi:\PP^1\rightarrow \PP^1$ of degree $4$ whose preperiodic
points are precisely the $x$-coordinates of the torsion
points of $E$.  Thus, Merel's Theorem would follow as a simple
corollary of the Morton and Silverman Conjecture for $N=1$
and $d=4$.
More generally, Fakhruddin has shown \cite{Fak2} that
the full Morton and Silverman Conjecture for $D=1$
would imply uniform boundedness of torsion for
abelian varieties.

The Conjecture seems to be very far from a proof.
However, there is growing evidence that it is valid,
at least in the simplest case, that $K=\QQ$, $N=1$, and
$\phi$ is a polynomial of degree $2$.
(The problem then reduces
to considering $\phi_c(z)=z^2 +c$, with $c\in\QQ$.)
In particular, the computations in \cite{Mor} and \cite{FPS}
show that $\phi_c$ never has a rational point of period $4$
or $5$ (respectively).
Moreover, Poonen showed in 1998
that if $\phi_c$ never has rational periodic points
of period greater than $5$, then it never has more than
$9$ rational preperiodic points \cite{Poo}.  (That bound, if
true, would be sharp, in light of the $c=-29/16$ example above.)
Those results all considered moduli spaces,
for fixed $n>m\geq 0$, of pairs $(c,x)$
such that $\phi_c^n(x)=\phi_c^m(x)$,
giving curves analogous to modular curves,
but with no known structure to take the place of a Hecke ring.  Instead,
the theorems were proven by delicate {\em ad hoc} computations
on the particular curves that arose.

Other researchers have found bounds for the longest possible period
of a $K$-rational periodic point by analyzing at a
prime of ``good reduction'' (see Definition~\ref{def:reduc} below);
see, for example, \cite{Fak1,MS1,Nar,Pez1,Pez2,Zieve}.
If $s$ is the total
number of primes of ``bad reduction,'' then
these results lead to bounds
on the order of at least $d^{s^{4D}}$ for the
number of $K$-rational periodic points (cf.~Corollary~B of
\cite{MS1}, for example).
Moreover, these results say virtually nothing about the more
general preperiodic points.

A different strategy (for the family $\phi_c(z)=z^2 + c$
for $c\in \QQ$) appeared in a 1997 theorem
of Call and Goldstine \cite{CaGol}.  They showed
that $\phi_c$ has at most $1+2^{s+2}$ rational preperiodic points,
where $s$ is the number of primes of bad reduction.
In their argument, they analyzed the dynamics
at the primes of bad, not good, reduction.
In particular, they studied the filled Julia set $\calK_v$
(see Definition~\ref{def:julia} below) at a bad prime $v$.
All preperiodic points except $\infty$ sit inside $\calK_v$,
which in turn lies in a union
of two $v$-adic disks, each of volume $1$.  For good $v$,
a single disk of volume $1$ suffices.  (A slightly different
condition holds at $2$ and at the archimedean (i.e., standard)
absolute value, of course.)  The bound
of $O(2^s)$ then follows naturally by an adelic argument.

In this paper, we will use a more detailed analysis of the filled
Julia sets to improve the above bounds substantially.
We will work only with polynomials and only in dimension $1$,
but we will allow arbitrary degree $d\geq 2$ and we will
work over an arbitrary global field $K$.
Following Call and Goldstine, we will study
dynamics over a global field $K$
by considering each of the associated
complete valued fields $\Cv$.  Thus, we
will end up studying both complex dynamics (if $\Cv\cong\CC$)
and the newer realm of $p$-adic and non-archimedean dynamics.
We refer the reader to \cite{Bea,CG,Mil} for expositions
on complex dynamics, and to \cite{Ben4,Ben5,Bez,Hs2,Riv1,Riv2,TVW}
for papers exploring various aspects of $p$-adic dynamics.

Our main result is the following.

\begin{mthm}
Let $K$ be a global field,
let $\phi\in K[z]$ be a polynomial of degree $d\geq 2$,
and let $s$ be the number of bad primes of $\phi$ in $K$.
Then the number of preperiodic points  of $\phi$
in $\PK$ is at most $O(s\log s)$.
\end{mthm}

A more precise statement appears in Theorem~\ref{thm:global}.
The big-$O$ constant is essentially $(d^2-2d+2)/\log d$
for large $s$.

The idea of the proof is to consider, for a given
prime $v$, the product
$$P_v = \prod_{i\neq j} |x_i - x_j|_v,$$
where $\{x_1,\ldots,x_N\}$ are finite
$K$-rational preperiodic points of $\phi$.
(The product $P_v$ is related to transfinite
diameters and capacities, as discussed in
Section~\ref{sect:cap1}.)
If $r_v$ is the diameter of the filled Julia set $\calK_v$,
then $P_v\leq r_v^{N(N-1)}$ naively; but in Lemmas~\ref{lem:EFbd}.a
and~\ref{lem:capbd}, we obtain
$P_v\leq r_v^{(d-1)N\log_d N}$ (with some correction
factors for $v$ archimedean).

The key, however, is our treatment of the prime $w$ with
filled Julia set $\calK_w$ of the largest diameter.
We partition $\calK_w$ into two pieces, and we show
in Lemmas~\ref{lem:EFbd}.b, \ref{lem:napart},
and~\ref{lem:apart}
that the corresponding product $P_w$ on each piece
satisfies $P_w\leq r_w^{(d-1)N(\log_d N - AN + B)}$
for certain simple constants $A$ and $B$.  The product $P$
of all the $P_v$'s, restricted to preperiodic points
in the given piece of $\calK_w$, is then bounded by
$r_w^{(d-1)NE}$, where $E=s\log_d N - AN + B$.  For $N$
slightly larger than $(1/A)s\log_d s$ (see
Lemma~\ref{lem:slogs}), we get $E<0$, so that $P<1$,
which contradicts the product formula for
the global field $K$.  Thus, we get a bound of about
$(1/A)s\log_d s$ on each piece; summing the two bounds
gives the Theorem.

Of course, the details are complicated.
In Sections~\ref{sect:nota} and~\ref{sect:jul}, we
will set terminology and recall fundamental
facts concerning local and global fields, bad primes, and
filled Julia sets.
In Section~\ref{sect:elem}, we will introduce notation
for certain expressions that will arise later,
and we will bound these expressions
in a series of technical but completely elementary Lemmas.
In Section~\ref{sect:cap1}, we will discuss transfinite diameters
and prove our first nontrivial bound for $P_v$,
for general bad primes.
In Sections~\ref{sect:cap2}
and~\ref{sect:cap3}, we will describe the partition
of the filled Julia set at a bad prime.
Finally, in Section~\ref{sect:global}, we will
state Theorem~\ref{thm:global}
and combine all the results from the preceding sections to prove it.

The author would like to thank Laura DeMarco, Andrew Granville,
Bjorn Poonen,
and Daniel Velleman for a number of helpful conversations relating
to details of certain auxiliary bounds and to some of the
references in the literature.  Many thanks also to
Matthew Baker for suggesting an improvement to the
archimedean case of Lemma~\ref{lem:capbd}, and for other
comments and stimulating discussions.

\section{Global Fields and Local Fields}
\label{sect:nota}

In this section we present the necessary fundamentals
from the theory of local and global fields.  We also
set some notational conventions for this paper.
Although this material
is well known to number theorists, we present it for the
convenience of dynamicists.  
See Section~B.1 of \cite{HS} or Section~4.4 of \cite{RV}
for more details on global fields and sets of absolute values;
see \cite{Gou,Kob} for expositions concerning the local
fields $\Cv$.

\subsection{Global fields and absolute values}
Throughout this paper, $K$ will denote a global field.  That is,
$K$ is either a number field (i.e., a finite extension of $\QQ$)
or a function field over a finite field (i.e., a finite extension
of $\FF_p(T)$ for some prime $p$).  We will write $M_K$ for the
set of standard absolute values on $K$.  That is, $M_K$ consists
of functions $|\cdot|_v:K\rightarrow \RR$ satisfying
$|x|_v\geq 0$ (with equality if and only if $x=0$),
$|xy|_v = |x|_v |y|_v$, and
$|x+y|_v\leq |x|_v + |y|_v$, for all $x,y\in K$.
(We will frequently abuse notation and write $v\in M_K$ when
our meaning is clear.)
Moreover, the absolute values in $M_K$ are chosen to
satisfy a {\em product formula}, which is to say that for each
$v\in M_K$, there is an integer $n_v\geq 1$ such that for all
$x\in K^{\times}$,
\begin{equation}
\label{eq:pfmla}
\prod_{v\in M_K} |x|_v^{n_v} = 1.
\end{equation}
(Implicit in the product formula is the fact that for any $x\in K^{\times}$,
we have $|x|_v=1$ for all but finitely many $v\in M_K$.)

All but finitely many $v\in M_K$ satisfy the ultrametric
triangle inequality
$$|x+y|_v \leq \max\{|x|_v,|y|_v\}.$$
(Note that this means $|n|_v\leq 1$ for all integers $n\in\ZZ$.)
Such $v$ are called {\em non-archimedean} absolute values;
the finitely many exceptions are called {\em archimedean} absolute
values.  If $K$ is a function field, then all absolute values
are non-archimedean.  If $K$ is a number field, then there are
archimedean absolute values, each of which, when restricted to $\QQ$,
is the familiar absolute value $|\cdot|$,
commonly denoted $|\cdot|_{\infty}$.
In fact, if $K$ is a number field, then
\begin{equation}
\label{eq:sumnv}
\sum_{\substack{v\in M_K, \\  v \text{ archimedean} } } n_v
= [K\col\QQ \hspace{0.75pt}] .
\end{equation}
Meanwhile, the non-archimedean absolute values in $M_K$
correspond to prime ideals
of the ring of integers of $K$.  For this reason, we frequently
refer to the absolute values $v\in M_K$ as primes of $K$, even
when $v$ is archimedean.

If $v$ is non-archimedean, then $|K^{\times}|_v$ is a discrete
subset of $\RR$, and we say that $v$ is a
{\em discrete valuation} on $K$.
In that case, let $\eps\in (0,1)$ be the largest absolute value
less than $1$ attained in $|K^{\times}|_v$, and choose
$\pi_v\in K$ with $|\pi_v|_v=\eps$.  Then $\pi_v$ is called
a {\em uniformizer} of $K$ at $v$, and we have
$|K^{\times}|_v = \{\eps^m : m\in\ZZ\}$.  Moreover, if $K$
is a number field, then $|\pi_v|_v^{-n_v}=p^f$ for some
prime number $p\in\ZZ$ and some positive integer $f$,
and $|\cdot|_v$ restricted to $\QQ$ is the usual
$p$-adic absolute value on $\QQ$.
In this case, we say that $v$ {\em lies above} the prime number $p$.

\subsection{Local fields}
For each $v\in M_K$, we can form the completion $K_v$
(often called the local field at $v$) of $K$ with
respect to $|\cdot|_v$.
We write $\Cv$ for the completion of an algebraic
closure $\bar{K}_v$ of $K_v$.  (The absolute value $v$ extends in
a unique way to $\bar{K}_v$ and hence to $\Cv$.)
The field $\Cv$ is then a complete and algebraically closed field.
If $v$ is archimedean, then $K_v$ is isomorphic either to $\RR$
(in which case we call $v$ a real prime)
or to $\CC$ (in which case we call $v$ a complex prime), and $\Cv\cong\CC$.
We will henceforth avoid the notation $K_v$, as we will soon
introduce the notation $\calK_v$ to denote a completely different
object in Section~\ref{sect:jul}.

If $v$ is non-archimedean, then $\Cv$ is not locally compact,
but it has other convenient properties not shared by $\CC$.
In particular, the disk
$\ints_v = \{c\in\Cv : |c|_v \leq 1\}$ forms a ring,
called the {\em ring of integers}, which has a unique
maximal ideal
$\maxid_v = \{c\in\Cv : |c|_v < 1\}$.  The quotient
$k_v=\ints_v/\maxid_v$ is called the {\em residue field} of $\Cv$.
The natural reduction map from $\ints_v$ to $k_v$, sending
$a\in\ints$ to $\bar{a}=a+\maxid_v\in k_v$, will
be used to define good and bad reduction of a polynomial
in Definition~\ref{def:reduc} below; but after proving a few
simple Lemmas about good and bad reduction, we will
not need to refer to $\ints_v$, $\maxid_v$, or $k_v$ again.

\subsection{Disks}
Let $\Cv$ be a complete and algebraically closed field with
absolute value $v$.  Given $a\in\Cv$ and $r>0$, we write
$$\Dbar(a,r) = \{x\in\Cv : |x-a|_v \leq r\}
\quad\text{and}\quad
D(a,r) = \{x\in\Cv : |x-a|_v < r\}$$
for the closed and open disks, respectively, of radius $r$
centered at $a$.  Note our convention that all
disks have positive radius.

If $v$ is non-archimedean and $U\subseteq\Cv$ is a disk,
then the radius of $U$ is unique; it is the same as the diameter
of the set $U$ viewed as a metric space.  However, any point
$b\in U$ is a center.  That is, if $|b-a|_v\leq r$, then
$\Dbar(a,r)=\Dbar(b,r)$, and similarly for open disks.
It follows that two disks intersect if and only
if one contains the other.  In addition, all disks are both
open and closed as topological sets; however, open disks and
closed disks can still behave differently in other ways.

Still assuming that $v$ is non-archimedean, the set
$|\Cv^{\times}|_v$ of absolute values actually attained
by elements of $\Cv^{\times}$ is usually not all
of $(0,\infty)$.  As a result, if
$r\in (0,\infty)\setminus |\Cv^{\times}|_v$,
then $D(a,r)=\Dbar(a,r)$ for any $a\in\Cv$.  However,
if $r\in |\Cv^{\times}|_v$, then $D(a,r)\subsetneq \Dbar(a,r)$.

\section{Bad Reduction and Filled Julia Sets}
\label{sect:jul}

The following definition originally appeared in \cite{MS1}.
We have modified it slightly so that ``bad reduction'' now
means not potentially good, as opposed to not good.
\begin{defin}
\label{def:reduc}
Let $\Cv$ be a complete, algebraically closed
non-archimedean field with
absolute value $|\cdot|_v$, ring of integers
$\ints_v = \{c\in\Cv : |c|_v \leq 1\}$,
and residue field $k_v$.
Let $\phi(z) \in \Cv(z)$ be a rational function
with homogenous presentation
$$\phi\left( [x,y]\right) = [f(x,y),g(x,y)],$$
where $f,g\in\ints_v[x,y]$ are relatively prime homogeneous polynomials
of degree $d=\deg\phi$, and at least one coefficient of $f$ or $g$
has absolute value $1$.  We say that $\phi$ has {\em good reduction}
at $v$
if $\bar{f}$ and $\bar{g}$ have no common zeros in
$k_v \times k_v$ besides
$(x,y)=(0,0)$.
We say that $\phi$ has {\em potentially good reduction} at $v$
if there is some linear fractional transformation
$h\in\PGL(2,\Cv)$ such that $h^{-1}\circ\phi\circ h$ has good
reduction.
If $\phi$ does
not have potentially good reduction,
we say it has {\em bad reduction} at $v$.
\end{defin}
Naturally, for $f[x,y]=\sum_{i=0}^d a_i x^i y^{d-i}$,
the reduction $\bar{f}[x,y]$ in Definition~\ref{def:reduc}
means $\sum_{i=0}^d \bar{a}_i x^i y^{d-i}$.
By convention, if $\Cv\cong\CC$ is archimedean,
we declare all rational functions in $\Cv(z)$ to have
bad reduction.

In this paper, we will consider only polynomial functions $\phi$
of degree at least $2$; that is, $\phi(z)=a_d z^d + \cdots + a_0$,
where $d\geq 2$, $a_i\in\Cv$, and $a_d\neq 0$.
If $\Cv$ is non-archimedean, then,
it is easy to check that $\phi$ has good reduction
if and only if $|a_i|_v\leq 1$ for all $i$ and $|a_d|_v = 1$.
In particular, by the product formula,
if $\phi\in K[z]$ for a global field $K$,
then there can be only finitely many primes $v\in M_K$ at which
$\phi$ has bad reduction.

The main focus of our investigation will be filled Julia
sets, which are standard objects of study in complex dynamics.
The motivating idea is that for a polynomial $\phi$, any point
$x$ of large enough absolute value will be sucked out to the
attracting fixed point at $\infty$ under iteration; thus, all
of the interesting dynamics involves points that do not escape
to $\infty$ under iteration.
Since we will be interested in both archimedean and non-archimedean
fields, we state the definition here more generally.
\begin{defin}
\label{def:julia}
Let $\Cv$ be a complete, algebraically closed
field with absolute value $|\cdot|_v$, and let
$\phi(z) \in \Cv[z]$ be a polynomial of degree $d\geq 2$.
The {\em filled Julia set} of $\phi$ at $v$ is
$$\calK_v = \{x\in\Cv: \{|\phi^n(x)|_v\}_{n\geq 1} \text{ is bounded} \}.$$
\end{defin}
We use the notation $\calK_v$ rather than $\calK_{\phi}$ because
in this paper, we will consider the polynomial $\phi\in K[z]$
to be fixed,
and we will study its filled Julia sets at various different
primes $v$ of $K$.

Non-archimedean filled Julia sets have been studied
in Section~5 of \cite{Ben5}, for example.  It is important to
note that while complex filled Julia sets are always compact,
their non-archimedean counterparts are not usually compact.
Fortunately, this technicality will not be an obstacle for
our investigations.

We note four fundamental properties of filled Julia sets.
First, $\calK_v$ is invariant
under $\phi$; that is, $\phi^{-1}(\calK_v) = \phi(\calK_v)=\calK_v$.
Second, all the finite preperiodic points of $\phi$
(that is, all the preperiodic points in $\PCv$ other than the
fixed point at $\infty$) are contained in $\calK_v$.
Third, if $U_0$ is a disk containing $\calK_v$, then
$$\calK_v = \bigcap_{n\geq 0} \phi^{-n}(U_0).$$
Finally, if the polynomial $\phi\in\Cv[z]$ has good reduction,
then $\calK_v=\Dbar(0,1)$.

Filled Julia sets have been studied extensively
in the archimedean case $\Cv=\CC$.  If $\phi_d(z)=z^d$,
then the (complex) filled Julia $\calK$ of $\phi_d$
is simply the closed unit disk $\Dbar(0,1)$.
Meanwhile, since the degree $d$
Chebyshev polynomial $\psi_d$ satisfies $\psi\circ h = h\circ \phi$,
where $h(z) = z + 1/z$, it follows
that the complex filled Julia set of $\psi_d$ is the interval
$[-2,2]$ in the real line.  These two examples are misleadingly
simple, however; most filled Julia sets are complicated fractal
sets.  For example, it is well known that
for $|c|>2$, the filled Julia set of
$\phi(z)=z^2 +c$ is homeomorphic to the Cantor set.
For many more complex examples (often in the form of the Julia set, which
is the boundary of the filled Julia set), see \cite{Bea,CG,Mil}.

There are not as many examples of non-archimedean filled Julia sets
in the literature.  For the convenience of the reader, we present
a few here.  More examples may be found in \cite{Ben5,Riv1}.

\begin{example}
\label{ex:cantor}
Given $\Cv$ non-archimedean and $d\geq 2$, fix $c\in\Cv$, and
consider $\phi(z)=z^d - c^{d-1}z$.  Assume for convenience
that $|d-1|_v = 1$.
If $|c|_v\leq 1$, then $\phi$ has good reduction,
and hence $\calK_v = \Dbar(0,1)$.  Thus, we consider $|c|_v > 1$;
let $r=|c|_v$ and $U=\Dbar(0,r)$.
Note that for $|x|_v>r$, we have $|\phi(x)|_v=|x|_v^d$,
so that $\phi^n(x)\rightarrow\infty$.
That is, $\calK_v\subseteq U_0$.

The set $\phi^{-1}(0)$ consists of $0$ and $d-1$ other points,
all distance $r$ from one another.  Using standard mapping properties
of non-archimedean polynomials (see, for example, Section~2
of \cite{Ben7}), it is not hard to show that $\phi^{-1}(U_0)$
consists of $d$ disks of radius $r^{2-d}$, each centered at one
of the points of $\phi^{-1}(0)$.  Moreover, each of these smaller
disks maps one-to-one and onto $U_0$, and in fact $\phi$ multiplies
distances by a factor of $r^{d-1} = |\phi'(0)|_v$
on each smaller disk.  It follows that $U_n=\phi^{-n}(U_0)$ is a union
of $d^n$ disks, each of radius $r^{1-(d-1)n}$.
(The sets $U_n$ are nested
so that each disk of $U_n$ contains exactly $d$ disks of $U_{n+1}$,
arranged so that any two are the maximal distance
$r^{1-(d-1)n}$ apart.)  Since the radii of the disks
approach zero, it is easy to verify that $\calK_v = \bigcap U_n$
is homeomorphic to a Cantor set on $d$ intervals.
\end{example}

\begin{example}
\label{ex:hsia}
Given $\Cv$ non-archimedean and $d\geq 3$, fix $a\in\Cv$, and
consider $\phi(z)=z^d - a z^{d-1}$.  Again, assume for convenience
that $|d-1|_v = 1$.
If $|a|_v\leq 1$, then $\phi$ has good reduction,
and hence $\calK_v = \Dbar(0,1)$.  As in the previous example,
then, we consider $|a|_v > 1$.  Let $r=|a|_v$;
once again, we have
$\calK_v\subseteq U_0=\Dbar(0,r)$.

This time, however, $\phi^{-1}(U_0)$ consists of only two
disks.  One, $W_1=\Dbar(a,r^{-(d-1)})$, is small
and maps one-to-one onto $U_0$; but the other, $W_2=\Dbar(0,1)$,
is comparatively large, and it maps $(d-1)$-to-$1$ onto $U_0$.
Because the disk $V=\Dbar(0,r^{-1/(d-2)})$ maps $(d-1)$-to-1 onto
itself, we can understand the shape of $\calK_v$ reasonably well.
In particular, $\phi^{-2}(U_0)$ consists of $d+2$ disks: two inside
$W_1$ (one is a preimage of $W_1$, and the other is a preimage
of $W_2$), and $d$ inside $W_2$ ($d-1$ are preimages of $W_1$, and
the last is a $(d-1)$-to-$1$ preimage of $W_2$).  If we continue
to take preimages, then each $U_n=\phi^{-n}(U_0)$ will be a union
of disks.  To compute $U_{n+1}$ from $U_n$, observe that
each disk of $U_n$ will
will have one preimage inside $W_1$ and (with one exception)
$d-1$ preimages inside $W_2$.  The one exception is the unique
disk of $U_n$ containing $0$; it has only one preimage
inside $W_2$, and that preimage maps $(d-1)$-to-one onto it.
Ultimately $\calK_v=\bigcap U_n$ will consist of the disk $V$
and all of its (infinitely many) preimages together with a
vaguely Cantor-like set at which the preimages of $V$ accumulate.

Thus, in contrast with Example~\ref{ex:cantor},
$\calK_v$ is neither a disk nor compact.  In general,
the filled Julia set of a polynomial of bad reduction over $\Cv$
will look something like this one.
However, the dynamics can be even more
complicated when there are regions on which
$\phi$ maps $n$-to-$1$
for some integer $n$ divisible by $p$, the characteristic
of the residue field $k_v$.
\end{example}

The preceding comments and examples made frequent reference to
disks $U_0$ containing $\calK_v$.  The smallest such disk will
be of particular importance to us.
The following Lemma shows the existence of the smallest disk
and gives a partial characterization of it.

\begin{lemma}
\label{lem:goodrad}
Let $\Cv$ be a complete, algebraically closed field with
absolute value $|\cdot|_v$.
Let $\phi\in\Cv[z]$ be a
polynomial of degree $d\geq 2$ with lead coefficient
$a_d\in\Cv$.  Denote by $\calK_v$
the filled Julia set of $\phi$ in $\Cv$.
Then:
\begin{list}{\rm \alph{bean}.}{\usecounter{bean}}
\item There is a unique smallest disk $U_0\subseteq\Cv$
  which contains $\calK_v$.
\item $U_0$ is a closed disk of some radius
  $r'_v \in |\Cv^{\times}|_v$, with $r'_v \geq |a_d|_v^{-1/(d-1)}$.
\item If $|\cdot|_v$ is non-archimedean, then
  $\phi$ has potentially good reduction if and only
  if $r'_v = |a_d|_v^{-1/(d-1)}$.
\end{list}
\end{lemma}

\begin{proof}
Choose $\alpha\in\Cv$ such that $\alpha^{d-1}=a_d$,
and let $\psi(z) = \alpha\phi(\alpha^{-1} z)$,
which is a monic
polynomial with filled Julia set $\alpha\calK_v$.
Given the scaling factors
of $|a_d|_v^{-1/(d-1)}=|\alpha|_v^{-1}$
in parts~(b) and~(c),
we may assume without
loss that $a_d=1$.

If $\Cv$ is archimedean, then $\Cv\cong \CC$.  In that case,
it is well known that $\calK_v$ is a compact set in the plane.
(See, for example, Lemma~9.4 of \cite{Mil}.)
Since $\calK_v$ is bounded, there is a unique
smallest disk containing $\calK_v$.
(See, for example, Exercise~3 in Appendix~I of \cite{YB}.)
Moreover, because $\calK_v$ is compact, this
disk must be closed, and so we denote it $\Dbar(b,r_v)$.

It is well known that the filled Julia set of a monic polynomial
over $\CC$ has capacity $1$; see, for example, Theorem~4.1
of \cite{BH}.  If $r_v<1$, then $\calK_v$ would fit inside
a disk of radius strictly smaller than $1$.  However, the capacity
of a disk in the plane is exactly its radius; as a result,
the capacity of $\calK_v$
would be strictly smaller than $1$, which is
a contradiction.  Therefore, $r_v\geq 1$, proving the
Lemma in the archimedean case.
(See Remark~\ref{rem:rouche} below for an alternate proof
not using capacity theory.)

If $\Cv$ is non-archimedean, then let $b\in\Cv$ be a fixed point
of $\phi$.  (Such $b$ exists because $\phi(z)-z$ is a polynomial
of degree $d\geq 2$.)  Clearly $b\in\calK_v$.  By
the coordinate change $z\mapsto z+b$, we may assume
that $b=0$.
Write $\phi(z) = z^d + a_{d-1}z^{d-1} + \cdots + a_1 z$.
Let $\tilde{r}_v=\max\{|a_i|_v^{1/(d-i)} : i=1,\ldots,d-1\}$
and $r_v = \max\{\tilde{r}_v, 1\}$;
note that $r_v\in |\Cv^{\times}|_v$.

If $r_v= 1$, then $\phi$ is a monic polynomial with
coefficients in $\ints_v$.
Hence, $\phi$ has good reduction;
with $U_0=\Dbar(0,1)$, the Lemma follows.

If $r_v>1$, then
the Newton polygon 
(see Section~6.5 of \cite{Gou} or Section~IV.3 of \cite{Kob})
for the equation $\phi(z)=0$ shows
that there is some $c\in\Cv$ with $|c|_v=r_v$
and $\phi(c)=0$.  In particular, any disk containing
$\calK_v$ must contain $\Dbar(0,r_v)$.

Moreover, if $|z|_v>r_v$, then
the $z^d$ term has larger absolute value than any other
term of $\phi(z)$, so that $|\phi(z)|_v = |z|_v^d$.
By induction, $|\phi^n(z)|_v = |z|_v^{d^n}$ for all $n\geq 1$.
It follows that $\calK_v\subseteq\Dbar(0,r_v)$.
By the previous paragraph, $\Dbar(0,r_v)$ is the smallest
disk containing $\calK_v$.

It only remains to show that if $r_v>1$, then $\phi$ cannot
have potentially good reduction.  However, after our coordinate
changes, $\phi$ is a monic polynomial without constant term.
The assumption that $r_v>1$ means that $|a_i|_v>1$ for some
coefficient $a_i$ of $\phi$.  Thus, by
Corollary~4.6 of \cite{Ben4}, $\phi$ cannot have good
reduction even after a change of coordinates.
\end{proof}

\begin{remark}
\label{rem:rouche}
The fact that $r_v\geq 1$ in the archimedean case can
also be proven directly, without reference to the power
of capacity theory.  The following alternate
argument was suggested to the author by Laura DeMarco.

Suppose that $\calK_v\subseteq \Dbar(a,r_v)$ for some $r_v<1$;
let $s=(1+r_v)/2$.  Since $\phi$ is a (monic) polynomial of degree
at least $2$, there is some radius $R>1$ such that for
all $z\in\CC$ with $|z-a|>R$, we have $|\phi(z)-a|>R$.

Let $A$ denote the annulus $\{z\in\CC: s\leq |z-a|\leq R\}$.
Every point of $A$
is attracted to $\infty$ under iteration of $\phi$.
Since $A$ is compact, there is some $n\geq 1$
such that $f(z) = \phi^n(z) - a$ has $|f(z)|>1$
for all $z\in A$.  Note that all $d^n$ zeros of $f$ lie
in $D(a,s)$.

Let $g(z) = (z-a)^{d^n} - f(z)$, which is a polynomial
of degree strictly less than $d^n$.  However,
for all $z\in\CC$ with $|z-a|=s$, we have
$$|f(z) + g(z)| = |z-a|^{d^n} = s^{d^n} < 1 < |f(z)|.$$
By Rouch\'{e}'s Theorem (noting that $g(z)\neq 0$
for $|z-a|=s$),
$f$ and $g$ have the same number of zeros in $D(a,s)$,
counting multiplicity.
That is a contradiction; thus, $r_v\geq 1$.
\end{remark}

The next two Lemmas give slightly more detailed information
about the filled Julia set for a polynomial of bad reduction over
a non-archimedean field.

\begin{lemma}
\label{lem:manydisk}
With notation as in Lemma~\ref{lem:goodrad},
suppose that
$\Cv$ is non-archimedean and
$r'_v>|a_d|_v^{-1/(d-1)}$.
Then $\phi^{-1}(U_0)$
is a disjoint union of closed disks $D_1,\ldots,D_{\ell}\subseteq U_0$,
where $2\leq \ell\leq d$.  Moreover, there are positive
integers $d_1,\ldots, d_{\ell}$ 
with $d_1+\cdots + d_{\ell} = d$ such that for each $i=1,\ldots, \ell$,
$\phi$ maps $D_i$ $d_i$-to-$1$ onto $U_0$.  That is, $\phi(D_i)=U_0$,
and every point $U_0$ has exactly $d_i$ pre-images in $D_i$, counting
multiplicity.
\end{lemma}

\begin{proof}
As in the previous proof, we may assume that
$\phi(z) = z^d + a_{d-1} z^{d-1} + \cdots + a_1 z$, that
$r'_v=\max\{|a_i|_v^{1/(d-i)}: 1\leq i \leq d-1\}>1$, and that
$U_0=\Dbar(0,r'_v)$.  As observed in that proof, $|\phi(z)|_v > |z|_v$
for $|z|_v > r'_v$.  Thus, $\phi^{-1}(U_0)\subseteq U_0$.
Moreover, $\phi(U_0)\supsetneq U_0$.

We now construct
the disks $D_1,\ldots, D_{\ell}$ inductively.
For each $i=1,2,\ldots$, suppose we already have $D_1,\ldots,D_{i-1}$,
and choose $b_i\in\phi^{-1}(U_0)\setminus(D_1\cup\cdots \cup D_{i-1})$.
(If this
is not possible, then skip to the next paragraph.)  By Lemmas~2.3 and~2.6
of \cite{Ben7}, there is a unique disk $D_i$ containing $b_i$ which
maps onto $U_0$, and this disk must be closed.  Since $D_j$ was
also unique for each $j<i$, the new disk $D_i$ must be disjoint
from $D_j$.  In addition, by Lemma~2.2 of \cite{Ben7}, $\phi$ maps
$D_i$ $d_i$-to-$1$ onto $U_0$, for some integer $d_i\geq 1$.

This process must stop with $\ell\leq d$,
because $\phi^{-1}(0)$ consists of exactly
$d$ points, counting multiplicity, and since each $d_i\geq 1$, at
least one must be contained in each $D_i$.
In fact, counting pre-images of $0$ also shows that
$d_1+\cdots+d_{\ell}=d$.

Finally, suppose that $\ell=1$; that is,
$\phi^{-1}(U_0)$ is a single disk
$D_1\subsetneq U_0$.  However, $\calK_v = \phi^{-1}(\calK_v)\subseteq D_1$,
contradicting the assumption
that $U_0$ is the smallest disk containing $\calK_v$.
Thus, we must have $\ell\geq 2$.
\end{proof}

%\begin{remark}
%The existence of the disks $D_i$, with degrees summing to $d$,
%can also be proven using
%well-known results from the theory of rigid analysis.  Specifically,
%the inverse image of a connected affinoid $U$ under a meromorphic function
%is a finite union of affinoids, each of which maps $d_i$-to-$1$
%onto $U$.  However, since we are dealing with a polynomial, we can
%prove the desired result with less machinery.
%\end{remark}

\begin{lemma}
\label{lem:minrad}
Let $K$ be a field with a discrete valuation $v$, and let
$\pi_v\in K$ be a uniformizer at $v$.  Let 
$\Cv$ be the completion of an algebraic closure of $K$.

Let $\phi(z) = a_d z^d + \cdots +a_0 \in K[z]$
be a polynomial of degree $d\geq 2$.
Denote by $\calK_v$ the filled Julia set of $\phi$ in $\Cv$,
and let $r'_v>0$ be the radius of the smallest disk in $\Cv$
containing $\calK_v$.  Suppose that
$r'_v> |a_d|_v^{-1/(d-1)}$.

If $\calK_v\cap K\neq \emptyset$, then
$$
|a_d|_v^{1/(d-1)} r'_v\geq
\begin{cases} 
|\pi_v|_v^{-1} & \text{ if } d=2,
\\
|\pi_v|_v^{-1/[(d-1)(d-2)]}
& \text{ if } d\geq 3.
\end{cases} 
$$
\end{lemma}

\begin{proof}
Given $b\in \calK_v \cap K$, we may replace $\phi$ by
$\phi(z+b)-b\in K[z]$, which is a polynomial with the
same degree and lead coefficient as $\phi$, but with
filled Julia set translated by $-b$.  In particular,
the radius $r'_v$ is preserved; so we may assume without
loss that $0\in\calK_v$.

As in the proof of Lemma~\ref{lem:goodrad},
choose $\alpha\in\Cv$ such that $\alpha^{d-1}=a_d$,
and let
$$\psi(z) = \alpha\phi(\alpha^{-1} z)
= z^d + \sum_{i=0}^{d-1} \alpha^{1-i} a_i z^i.$$
Then $\psi$ is a monic
polynomial with filled Julia set $\calK'_v=\alpha\calK_v$;
however, $\psi$ may not be defined over $K$.
Still, the radius $r_v$ of the smallest disk containing
$\calK'_v$ satisfies $r_v>1$, by hypothesis.

Let $j$ be the largest
index between $0$ and $d-1$ that maximizes
$\lambda_j=|\alpha^{1-j} a_j|_v^{1/(d-j)}$.
Note that $\lambda_j > 1$;
for if not, then $\psi$
has good reduction, contradicting Lemma~\ref{lem:goodrad}.c.
The Newton polygon for the equation $\psi(z)=0$
shows that there is some $\beta\in\Cv$ with
$\psi(\beta)=0$ and $|\beta|_v = \lambda_j$.
We have $0,\beta\in\calK'_v$;
hence,
$r_v\geq \lambda_j$.

If $j=0$, then a simple induction shows that
$|\psi^n(0)|_v = |\alpha a_0|_v^{d^{n-1}}$ for $n\geq 1$.
Since $|\alpha a_0|_v>1$, this contradicts the
hypothesis that $0\in\calK'_v$.

Thus, $1\leq j\leq d-1$, and we write $|a_d|_v = |\pi_v|_v^{e_1}$
and $|a_j|_v = |\pi_v|_v^{e_2}$; note that $e_1,e_2\in\ZZ$.
Our assumptions say that
$$r_v\geq \lambda_j
= |\alpha^{1-j} a_j|_v^{1/(d-j)}
= |\pi_v|_v^f > 1,
\qquad\text{where}\qquad
f = \frac{1}{d-j}\left( \frac{(1-j)e_1}{d-1} + e_2\right) < 0.$$
If $j=1$, then $f=e_2/(d-1)\leq -1/(d-1)$, which proves the Lemma
for $d=2$.
If $2\leq j \leq d-1$, then
$f \leq -1/[(d-1)(d-j)]\leq -1/[(d-1)(d-2)]$, and we are done.
\end{proof}

\begin{remark}
The bounds of Lemma~\ref{lem:minrad} are sharp.  Indeed,
if $d=2$, then the polynomial $\phi(z)=z^2 - \pi_v^{-1} z$
has $0,\pi_v^{-1}\subseteq \calK_v$.  Because $\phi^n(z)\rightarrow\infty$
for $|z|_v > |\pi_v|_v^{-1}$, the smallest disk containing $\calK_v$
is $\Dbar(0,|\pi_v|_v^{-1})$, so that
$r_v = |\pi_v|_v^{-1}$.

Similarly, if $d\geq 3$, then the polynomial
$\phi(z)=\pi_v^d z^d - \pi_v z^2$ has $0\in\calK_v$.
Choosing a $(d-1)$-st root of $\pi_v$, we get
an associated monic conjugate 
$\psi(z)=z^d - \pi_v^{-1/(d-1)} z^2$, from which it
is easy to compute
$|\pi_v^d|_v^{1/(d-1)}r'_v= r_v = |\pi_v|_v^{-1/[(d-1)(d-2)]}$.
\end{remark}

\section{Elementary Computations}
\label{sect:elem}

We will now define and bound certain
integer quantities that will appear as exponents
in the rest of the paper.  
The reader is encouraged to read the
statements of Definition~\ref{def:E},
Lemma~\ref{lem:EFbd}, and Lemma~\ref{lem:slogs}
but to skip the proofs, which are 
tedious but completely elementary,
until after seeing their use in Theorem~\ref{thm:global}.

We will write $\log_d x$ to denote the logarithm of $x$
to base~$d$.

\begin{defin}
\label{def:E}
Let $N\geq 0$ and $d\geq 2$ be integers.  We define
$E(N,d)$ to be twice the sum of all base-$d$ coefficients of all integers
from $0$ to $N-1$.  That is,
$$E(N,d) = 2\sum_{j=0}^{N-1} e(j,d),
\qquad\text{where}\qquad
e\left(\sum_{i=0}^M c_i d^i , d\right) = \sum_{i=0}^M c_i,$$
for $c_i\in\{0,1,\ldots, d-1\}$.

Moreover, if $m$ is an integer satisfying
$1\leq m\leq d$,
we may write $N=c_0 + m k$ for unique integers $c_0\in\{0,1,\ldots, m-1\}$
and $k\geq 0$.  We then define
$$e(N,m,d) = c_0 + e(k,d) -(d-m)k
\qquad\text{and}\qquad
f(N,m,d) = c_0 + e(k,d),$$
and
$$E(N,m,d) = 2\sum_{j=0}^{N-1} e(j,m,d)
\qquad\text{and}\qquad
F(N,m,d) = 2\sum_{j=0}^{N-1} f(j,m,d).$$
\end{defin}

We declare $E(N,d)=E(N,m,d)=F(N,m,d)=0$ for $N\leq 1$.
Clearly, $E(N,d)$ and $F(N,m,d)$ are always positive for $N\geq 1$;
but for $N$ large and $m<d$,
$E(N,m,d)$ is negative.
Note that $e(N,d,d)=f(N,d,d)=f(N,1,d)=e(N,d)$, and therefore
\begin{equation}
\label{eq:EFid}
E(N,d,d)=F(N,d,d)=F(N,1,d)=E(N,d).
\end{equation}

We will need 
the following two auxiliary Lemmas.

\begin{lemma}
\label{lem:EandF}
Let $N,m,d$ be integers satisfying
$N\geq 1$, $d\geq 2$, and $1\leq m\leq d$.  Write
$N=c+mk$ with $0\leq c\leq m-1$ and $k\geq 0$.
Then:
\begin{list}{\rm \alph{bean}.}{\usecounter{bean}}
\item
  $F(N,m,d) = (m-c)E(k,d) + cE(k+1,d) + (m-1)N - c(m-c)$.
\item
  $\dsps E(N,m,d) = F(N,m,d) - \frac{(d-m)}{m}[N^2 - mN + c(m-c)]$.
\item If $N\leq m$, then
  $E(N,m,d) = F(N,m,d) = N(N-1)$.
\end{list}
\end{lemma}

\begin{proof}
Writing an arbitrary integer $j\geq 0$ as $j=i+m\ell$ for
$0\leq i\leq m-1$, we compute
\begin{align*}
F(N,m,d) & = 2 \sum_{j=0}^{N-1}f(j,m,d)
= 2\sum_{i=0}^{c-1}\sum_{\ell=0}^k f(i+m\ell,m,d)
+ 2\sum_{i=c}^{m-1}\sum_{\ell=0}^{k-1} f(i+m\ell,m,d)
\\
& = 2 \sum_{i=0}^{c-1}\sum_{\ell=0}^k (i + e(\ell,d))
+ 2 \sum_{i=c}^{m-1}\sum_{\ell=0}^{k-1} (i + e(\ell,d))
\\
& = \sum_{i=0}^{c-1} \left[2(k+1)i + E(k+1,d)\right]
+ \sum_{i=c}^{m-1} \left[2ki + E(k,d)\right]
\\
& =c E(k+1,d) + (m-c) E(k,d) + (k+1)c(c-1) + km(m-1) - kc(c-1).
\end{align*}
Part (a) of the Lemma now follows by rewriting the last
three terms as
$$ c(c-1) + mk(m-1) = c(c-m) + (c+mk)(m-1) = (m-1)N -c(m-c).$$

Next, we compute
\begin{align*}
E(N,m,d) &= 2 \sum_{j=0}^{N-1}e(j,m,d)
= 2\sum_{j=0}^{N-1}f(j,m,d) - 2(d-m)\left[
m\sum_{\ell=0}^{k-1} \ell + \sum_{j=0}^{c-1} k \right]
\\
& = F(N,m,d) - k(d-m)[ m(k-1) + 2c ]
\end{align*}
Writing $k=(N-c)/m$, the last term becomes
$$-\frac{(N-c)}{m}(d-m)(N+c-m)
= - \frac{(d-m)}{m}[N^2-mN + c(m-c)],$$
proving part (b).
Finally, part (c) is immediate from the fact that
$e(j,m,d)=f(j,m,d)=j$ for $0\leq j \leq m-1$.
\end{proof}

\begin{lemma}
\label{lem:concave}
Let $N,m,d$ be integers satisfying
$N\geq 1$, $d\geq 2$, and $1\leq m\leq d$.  Write
$N=c+mk$ with $0\leq c\leq m-1$ and $k\geq 0$.
Then:
\begin{list}{\rm \alph{bean}.}{\usecounter{bean}}
\item
$\dsps (m-c) \log_d\left(\frac{mk}{N}\right)
+ c \log_d\left(\frac{mk+m}{N}\right) \leq 0$.
\item If $N\geq d$, then
$\dsps (d-1)\log_d\left(\frac{mk+m}{N}\right) - (m-c)\leq 0$.
\end{list}
\end{lemma}

\begin{proof}
The function $\log_d(x)$ is of course concave
down.  Letting $x_1=mk/N$ and $x_2=(mk+m)/N$, then, we have
$x_1\leq 1< x_2$, and therefore $\log_d(1)\geq L(1)$, where
$$L(x) = \frac{1}{x_2 - x_1}
\left[(x_2 - x) \log_d(x_1) + (x-x_1)\log_d(x_2) \right]$$
is the line through $(x_1,\log_d(x_1))$ and $(x_2,\log_d(x_2))$.
That is,
$$0 \geq \frac{1}{m}\left[(m-c)\log_d\left(\frac{mk}{N}\right)
+ c \log_d\left(\frac{mk+m}{N}\right)\right],$$
proving part (a).
For part (b), we have
$$(d-1)\log_d\left(\frac{mk + m}{N}\right) =
\frac{(d-1)}{\log d} \cdot \log\left( 1 + \frac{m-c}{N} \right)
\leq \frac{(d-1)}{\log d} \cdot \frac{(m-c)}{N}.$$
However, $\log d = -\log [1 - (d-1)/d ] \geq (d-1)/d$,
and since $N\geq d$,
\[
(d-1)\log_d\left(\frac{mk+m}{N}\right)
\leq (d-1) \cdot \frac{d}{d-1}\cdot \frac{m-c}{N}
=\frac{d}{N}(m-c) \leq (m-c).
\qedhere
\]
\end{proof}

\begin{lemma}
\label{lem:EFbd}
Let $N,m,d$ be integers satisfying $N\geq 1$, $d\geq 2$,
and $1\leq m\leq d-1$.
Then:
\begin{list}{\rm \alph{bean}.}{\usecounter{bean}}
\item
$E(N,d)\leq (d-1) N \log_d N $, with equality if $N$ is
a power of $d$.
\item
$\dsps E(N,m,d)\leq (d-1) N
\left[ \log_d N + 1 - \log_d m - \frac{(d-m)}{m(d-1)} N\right]$,
with equality if $N/m$ is a power of $d$.
\item
$F(N,m,d)\leq (d-1)N \log_d N$.
\item
For $N\geq m$,
$\dsps F(N,m,d)\leq (d-1) N
\left[ \log_d N - \log_d m + \frac{m-1}{d-1}\right]$,
with equality if $N/m$ is a power of $d$.
\end{list}
\end{lemma}

\begin{proof}  Fix $d\geq 2$.
If $N=1$, then both sides of part~(a) are clearly zero.
If $2\leq N\leq d$, then $E(N,d)=N(N-1)$ by Lemma~\ref{lem:EandF}.c
(with $m=d$) and equation~\eqref{eq:EFid}.
Because $(\log x)/(x-1)$
is a decreasing function of the real variable $x> 1$,
we have
$(\log d)/(d-1) \leq (\log N)/(N-1)$, with
equality for $N=d$.
Part~(a) then follows for $1\leq N\leq d$.

For $N\geq d+1$, we proceed by induction
on $N$, assuming part~(a) holds for all positive integers up to $N-1$.
Write $N=c+ dk$, where $0\leq c\leq d-1$, so that
$1\leq k \leq N-2$.
By Lemma~\ref{lem:EandF}.a (with $m=d$) and equation~\eqref{eq:EFid},
we have
\begin{align*}
E(N,d)
& =(d-c) E(k,d) + c E(k+1,d) + (d-1)N - c(d-c)
\\
& \leq (d-c)(d-1) k \log_d k + c(d-1)(k+1)\log_d(k+1)
%\\
%& \phantom{{}\leq (d-c)(d-1)k\log_d k }
+ (d-1)N - c(d-c)
\\
& = (d-c)(d-1)k\log_d(dk) + c(d-1)(k+1)\log_d(dk+d) - c(d-c)
\end{align*}
where the final equality is because
$N=(d-c)k + c(k+1)$, and
the inequality (which is equality if $N$ is a power of $d$)
is by the inductive hypothesis, since
$k,k+1\leq N-1$.
More generally, adding and subtracting $(d-1)N\log_d N$, we have
\begin{align*}
E(N,d) & \leq (d-1)N\log_d N+ 
(d-c)(d-1)k\log_d\left(\frac{dk}{N}\right)
\\
& \phantom{{} \leq(d-1)N\log_d N}
+ c(d-1)(k+1)\log_d\left(\frac{dk+d}{N}\right) - c(d-c)
\\
& = (d-1)N\log_d N
+ c\left[ (d-1)\log_d\left(\frac{dk+d}{N}\right) - (d-c) \right].
\\
& \phantom{{} \leq(d-1)N\log_d N}
+ (d-1)k \left[ (d-c) \log_d\left(\frac{dk}{N}\right)
+ c \log_d\left(\frac{dk+d}{N}\right)\right]
\end{align*}
By Lemma~\ref{lem:concave} with $m=d$,
the quantities in brackets are nonpositive,
and part~(a) follows.

If $m=1$, then parts~(c--d) are immediate from
part~(a) and equation~\eqref{eq:EFid}.  Moreover,
by Lemma~\ref{lem:EandF}.a--b (with $m=1$)
and part~(a),
$$E(N,1,d) = E(N,d) - (d-1) N (N-1)
\leq (d-1) N [ \log_d N + 1 - N],$$
with equality if $N$ is a power of $d$.
This is exactly part~(b) for $m=1$.
Thus, we may assume for the remainder of the proof that $2\leq m\leq d$.

We now turn to part~(d).  If $N=m$, then
by Lemma~\ref{lem:EandF}.c, we have $F(m,m,d)\leq m(m-1)$,
which exactly equals the desired right hand side.
For $N\geq m+1$, write $N=c+mk$, where $k\geq 1$ and
$0\leq c\leq m-1$.  By Lemma~\ref{lem:EandF}.a,
\begin{equation}
\label{eq:F1}
F(N,m,d) =
(m-c) E(k,d) + c E(k+1,d) + (m-1)N - c(m-c).
\end{equation}
If $m+1\leq N\leq d-1$, then $k\leq d-1$, so that
by Lemma~\ref{lem:EandF}.c, equation~\eqref{eq:F1}
becomes
\begin{align*}
F(N,m,d) &= (m-c)k(k-1) + ck(k+1) + (m-1)N - c(m-c)
\\
&= mk^2 -mk + 2ck + (m-1)N - c(m-c)
= (m+k-1)N - (m-c)(k+c)
\\
&= m^{-1} \left[ (N + m^2 - c -m)N - (m-c)(N-c + cm) \right]
\\
&= m^{-1} \left[ (N + m^2 -2m)N - c(m-c)(m-1)\right]
\leq N\left[(N/m) + m -2\right],
\end{align*}
where we have substituted $k=(N-c)/m$ along the way.
Thus, we must show
$$N\left[(N/m) + m -2\right]
\leq N\left[(d-1)\log_d(N/m) + m-1 \right].$$
Equivalently, we must show
$$ \frac{\log d}{d-1} \leq \frac{\log(N/m)}{(N/m)-1},$$
which is true because $(\log x)/(x-1)$ is a decreasing
function of $x> 1$, and $1< N/m < d$.

If $N\geq d$ in part~(d), then
$k\geq 1$. Applying part~(a) to equation~\eqref{eq:F1},
we obtain
\begin{multline}
\label{eq:F3}
F(N,m,d) \leq (m-c)(d-1) k \log_d k + c(d-1)(k+1)\log_d(k+1) \\
+ (m-1)N -c(m-c),
\end{multline}
with equality if $c=0$ and $k$ is a power of $d$, whence
we immediately obtain the statement of the Lemma
for $N=md^i$.  More generally,
\eqref{eq:F3} becomes
\begin{align*}
F(N,m,d) & \leq (d-1)N\left(\log_d \frac{N}{m} + \frac{m-1}{d-1}
\right)
+ (m-c)(d-1)k\log_d\left(\frac{mk}{N}\right)
\\
& \phantom{{} \leq (d-1)N\left(\log_d \frac{N}{m} \right.}
+ c(d-1)(k+1)\log_d\left(\frac{mk+m}{N}\right) - c(m-c)
\\
& = (d-1)N\left(\log_d \frac{N}{m} + \frac{m-1}{d-1} \right)
+ c\left[ (d-1)\log_d\left(\frac{mk+m}{N}\right) - (m-c) \right]
\\
& \phantom{{} \leq (d-1)N\left(\log_d \frac{N}{m} \right.}
+ (d-1)k \left[ (m-c) \log_d\left(\frac{mk}{N}\right)
+ c \log_d\left(\frac{mk+m}{N}\right)\right].
\end{align*}
Part~(d) now follows from Lemma~\ref{lem:concave}, as before.

For part~(c), if $1\leq N\leq m$, then by
part~(a) and Lemma~\ref{lem:EandF}.c,
$$F(N,m,d) =N(N-1) = E(N,d) \leq (d-1)N \log_d N,$$
as desired.
The remaining case, that $N\geq m$, will follow from part~(d)
provided
$$m-1\leq (d-1)\log_d m.$$
However, this is the same as showing that
$\log d/(d-1) \leq \log m/(m-1)$,
which once again follows from the fact that $(\log x)/(x-1)$ is
decreasing for $x> 1$.

Last, we turn to part~(b).
If $1\leq N\leq m$, then
$E(N,m,d) = N(N-1)$
by Lemma~\ref{lem:EandF}.c.
Thus, we wish to show that
$$N-1 \leq (d-1)\left( 1 + \log_d \frac{N}{m}\right) -(d-m) \frac{N}{m},$$
which is to say
$$\frac{dN}{m} - 1 \leq (d-1) \log_d \left( \frac{dN}{m}\right),$$
with equality when $N=m$.
Yet again, this inequality
follows immediately from the facts that $(\log x)/(x-1)$
is decreasing for $x> 1$ and that $1\leq dN/m\leq d$.

It only remains to consider $N\geq m+1$.
Writing $N=c+mk$, where $k\geq 1$ and $0\leq c\leq m-1$,
and invoking Lemma~\ref{lem:EandF}.b, we have
\begin{align*}
E(N,m,d) & = F(N,m,d) - \frac{(d-m)}{m}[N^2-mN + c(m-c)]
\\
& \leq F(N,m,d) -\frac{(d-m)}{m}N^2 + (d-m)N,
\end{align*}
with equality for $c=0$.
By part~(d), we obtain
\begin{align*}
E(N,m,d) & \leq 
(d-1) N \left[ \log_d N - \log_d m + \frac{m-1}{d-1}\right]
- \frac{(d-m)}{m} N^2 + (d-m)N
\\
& = (d-1) N \left[ \log_d N - \log_d m + 1  - \frac{d-m}{m(d-1)}N\right],
\end{align*}
with equality if $N$ is of the form $N=md^i$.
\end{proof}

Besides the preceding integer quantities and their bounds,
we will need the following bound involving a certain family
of real-valued functions.

\begin{lemma}
\label{lem:slogs}
Let $d\geq 2$ be an integer, and let
$A,B,t$ be positive real numbers such that
$$(d-1)A \geq d^{B-1}
\qquad \text{and} \qquad
t\geq 1.$$
Define $\eta:(0,\infty)\rightarrow \RR$
by
$$\eta(x) = t\log_d x - Ax + B.$$
Set the real number $M(A,B,t)$ to be
$$M(A,B,t) =
\frac{t}{A}\left(\log_d t + \log_d(\max\{1,\log_d t\}) + 3\right) .
$$
Then $\eta(x)<0$ for all $x\geq M(A,B,t)$.
\end{lemma}

\begin{proof}
By differentiating, we see that $\eta$ is decreasing
for $x\geq t/(A \log d)$, and hence for $x\geq M(A,B,t)$.
Thus, it suffices to show that $\eta(M(A,B,t))<0$.

First, suppose that $t<d$.  Then
\begin{align*}
\eta(M(A,B,t)) & = t\log_d t + t\log_d\left[ A^{-1}(\log_d t + 3)\right]
- t\log_d t - 3t + B
\\
& = t\log_d\left[ A^{-1} d^{B-3}(\log_d t + 3) \right]
-B(t-1)
\leq t\log_d\left[ (d-1)(\log_d t + 3)/d^2 \right] ,
\end{align*}
where the inequality is because $A^{-1}d^{B-1} \leq (d-1)$ and $B>0$, by
hypothesis.  Since $t<d$, the quantity inside square brackets
is strictly less than $4(d-1)/d^2\leq 1$.  Thus,
$\eta(M(A,B,t))< t\log_d(1) = 0$, and we are done.

Second, if $t\geq d$, then by a similar computation,
\begin{align*}
\eta(M(A,B,t)) &=
t\log_d\left[ A^{-1}d^{B-3}u^{-1}(u + \log_d u + 3) \right] -B(t-1)
\\
& < t\log_d\left[ (d-1)(u + \log_d u + 3)/(d^2 u) \right]
\end{align*}
where $u=\log_d t$.
Writing
$H(u) = (d-1)(u + \log_d u + 3)/(d^2 u)$,
it suffices to show that $H(u)\leq 1$ for $u\geq 1$.
Differentiating, it is easy to see that $H$ is decreasing
for such $u$.  Since $H(1)=4(d-2)/d^2 \leq 1$, we are done.
\end{proof}

\section{Transfinite Diameters and Bad Primes}
\label{sect:cap1}

Given a metric space $X$ and an integer $N\geq 2$,
the $N^{\text{th}}$ diameter of $X$ is defined to be
$${\mathbf d}_N(X)
= \sup_{x_1,\ldots,x_N \in X} \prod_{i\neq j}d_X(x_i,x_j)^{1/[N(N-1)]},$$
which measures the maximal average distance between any two of $N$ points
in $X$.
(See \cite{GJ}, for example, 
for a computation of the $N^{\text{th}}$ diameter
of the interval $[0,1]$.)
This quantity is usually used to define the
{\em transfinite diameter} of $X$,
$${\mathbf d}(X) = \lim_{N\rightarrow\infty} {\mathbf d}_N(X),$$
which converges because $\{ {\mathbf d}_N(X)\}_{N\geq 2}$ is a
decreasing sequence.  
If $X$ is a nice enough (e.g., compact) subset
of a valued field, then the transfinite
diameter coincides with the Chebyshev constant and the
logarithmic capacity of $X$; see Section~5.4 of \cite{Ami},
or Chapters~3 and~4 of \cite{Rum}.
Baker and Hsia used this equality in \cite{BH} to compute
the transfinite diameter
of filled Julia sets of polynomials, even when
those sets were not compact.  (Their result
of $|a_d|_v^{-1/(d-1)}$, where $d$ is the degree and $a_d$
the lead coefficient of the polynomial, was already well known
for $\Cv=\CC$.)  See \cite{Rum} for more on transfinite diameters
and capacities in $\Cv$.

However, in this paper we will be interested in the
$N^{\text{th}}$ diameters ${\mathbf d}_N(X)$ themselves,
rather than the transfinite diameter.  In particular, 
the following Lemma contains our main bound for
${\mathbf d}_N(\calK_v)^{N(N-1)}$, where $\calK_v$ is
the filled Julia set of a polynomial $\phi\in\Cv[z]$.
The proof uses an
estimate involving van der Monde determinants
similar to a bound that appears in the proof
of Lemme~5.4.2 in \cite{Ami}.

\begin{lemma}
\label{lem:capbd}
Let $\Cv$ be a complete, algebraically closed field with
absolute value $|\cdot|_v$.
Let $\phi\in\Cv[z]$ be a
polynomial of degree $d\geq 2$ with lead coefficient
$a_d\in\Cv$.  Denote by $\calK_v$
the filled Julia set of $\phi$ in $\Cv$, and let $r'_v$
be the radius of the smallest disk that contains $\calK_v$.
Set $r_v = |a_d|_v^{1/(d-1)} r'_v$.

Then for any integer $N\geq 2$ and any set $\{x_1,\ldots,x_N\}
\subseteq \calK_v$ of $N$ points in $\calK_v$,
$$\prod_{i\neq j} |x_i - x_j|_v \leq
|a_d|_v^{-N(N-1)/(d-1)} \max\{1, |N|_v^N\} r_v^{E(N,d)},$$
where $E(N,d)$ is twice
the sum of all base-$d$ coefficients of all integers
from $0$ to $N-1$, as in Definition~\ref{def:E}.
\end{lemma}

\begin{proof}
%We begin by reducing to the case that $\phi$ is monic.
Choose $\alpha\in\Cv$ such that $\alpha^{d-1}=a_d$,
and let $\psi(z) = \alpha\phi(\alpha^{-1} z)$.  Then $\psi$ is a monic
polynomial with filled Julia set
$\calK'_v=\alpha\calK_v$, and the smallest disk containing
$\calK'_v$ has radius $r_v$.  If the Lemma holds
for $\psi$, then for $x_1,\ldots,x_N\in\calK_v$, we have
$\alpha x_i \in \calK'_v$, and therefore
$$\prod_{i\neq j} |x_i - x_j|_v
= |\alpha|_v^{-N(N-1)} \prod_{i\neq j} |\alpha x_i - \alpha x_j|_v
 \leq |\alpha|_v^{-N(N-1)}
\max\{1, |N|_v^N\} r_v^{E(N,d)},$$
as desired.
Thus, it suffices to prove the Lemma in the case that $\phi$
is monic.

We will now construct a sequence $\{f_j\}_{j=1}^{\infty}$ of
monic polynomials over $\Cv$ such that each $f_j$ has
degree $j$ and such that $|f_j(x)|_v$ is not especially large
for any $x\in\calK_v$.

First, let $\Dbar(a,r_v)$ be the smallest disk
containing $\calK_v$, where $a\in\Cv$ and $r_v$
is as in the statement of the Lemma.
For any integer $j\geq 0$ written in base-$d$
notation as
$$j= c_0 + c_1 d + c_2 d^2 + \cdots + c_M d^M,$$
with $c_i\in\{0,1,\ldots, d-1\}$, define
$$f_j(z) = \prod_{i=0}^M [\phi^i(z) - a]^{c_i}.$$
Clearly, $f_j$ is monic of degree $j$.  Moreover, for $x\in\calK_v$,
we have $\phi^i(x)\in\calK_v$, and therefore
$$|f_j(x)|_v \leq \prod_{i=0}^M r_v^{c_i} = r_v^{e(j,d)},$$
where $e(j,d)$ is as in Definition~\ref{def:E}.

Given $x_1,\ldots, x_N\in \calK_v$, denote
by $V(x_1,\ldots,x_N)$ the corresponding van der Monde matrix
(i.e., the $N\times N$ matrix with
$(i,j)$ entry $x_i^{j-1}$).  Recall that
$$\prod_{i\neq j} |x_i - x_j|_v = |\det V(x_1,\ldots, x_N)|_v^2.$$
Because $f_{N-1}$ is monic,
we may replace the last column of the matrix
by a column with
entry $f_{N-1}(x_i)$ in the $i$th row, without changing
the determinant.
We may then replace the second to last column
by a column with entry $f_{N-2}(x_i)$ in the $i$th row,
and so on.  Thus, if we denote by $A(x_1,\ldots,x_N)$
the matrix with $(i,j)$ entry $f_{j-1}(x_i)$, then
$$\det V(x_1,\ldots, x_N) = \det A(x_1,\ldots, x_N).$$
If $\CC_v=\CC$ is archimedean, then by Hadamard's inequality
applied to the columns of $A$,
$$|\det A(x_1,\ldots, x_N)|^2
\leq \prod_{j=0}^{N-1} \left(|f_j(x_1)|^2 + \cdots + |f_j(x_N)|^2\right)
\leq \prod_{j=0}^{N-1} N r_v^{2e(j,d)}
= N^N r_v^{E(N,d)}.$$
Similarly, if $\CC_v$ is non-archimedean, then by
the non-archimedean version of Hadamard's inequality
(see, for example, \cite{Ami}, Preuve du Lemme~5.3.4),
we have
\[
|\det A(x_1,\ldots, x_N)|^2
\leq \prod_{j=0}^{N-1} \max_{i=1,\ldots,N} |f_j(x_i)|^2
\leq \prod_{j=0}^{N-1} r_v^{2e(j,d)}
= r_v^{E(N,d)}. \qedhere
\]
\end{proof}

\begin{remark}
We can recover the Baker and Hsia bound
${\mathbf d}(\calK_v) \leq |a_d|_v^{-1/(d-1)}$ immediately
from Lemmas~\ref{lem:capbd} and~\ref{lem:EFbd}.a.
(The opposite inequality is more subtle, however.)
\end{remark}

\begin{remark}
\label{rem:sharp1}
There are many cases for which the bound of Lemma~\ref{lem:capbd}
is sharp.  In particular, for non-archimedean $v$, degree
$d\geq 2$ with $|d-1|_v =1$, and $c\in \Cv$
with $|c|>1$, recall that the function $\phi(z)=z^d - c^{d-1}z$ of
Example~\ref{ex:cantor} has $\calK_v$ homeomorphic
to a Cantor set on $d$ pieces.  For arbitrary $N\geq 2$, one
can distribute $N$ points in $\calK_v$ in the following way.
Write $N=\sum_{i=0}^{M} c_i d^i$, and put $c_M$ points
in each of the $d^M$ pieces at level $M$, maximally far apart in
each piece; then put $c_{M-1}$ in each of the $d^{M-1}$ pieces
at level $M-1$, each as far as possible from the existing
points; and so on.  Keeping track of the radii of
the disks at each level, one can show that
$\prod_{i\neq j} |x_i-x_j|_v = r_v^{E(N,d)}$
exactly.

In many other cases, however, the bound is not quite sharp,
though it appears to be approximately the right order of
magnitude.  In the archimedean case, of course, the Hadamard
inequality introduces some error.  Still, the greater factor
seems to be the choice of the monic polynomial $f_j$.  When
$j$ is a power of $d$, computations suggest that our choice
of $f_j$ is very close to sharp, if not actually sharp.
However, when $j$ is not a power of $d$, our construction
of $f_j$ as a product of smaller factors is in general not
optimal, even in the non-archimedean setting.
For example, if $\phi(z) = z^3 - az^2$ is the
map of Example~\ref{ex:hsia} (non-archimedean, with $d=3$,
$|a|_v>1$,
and $|2|_v=1$), then the function $f_6(z)=(\phi(z))^2$
of the proof has $|f_6(z)|_v$ growing as large as $r^2$
on $\calK_v$;
but the function $\tilde{f}_6(z) = (\phi(z))\cdot(\phi(z)-a)$
has $|\tilde{f}_6(z)|_v\leq r$.  Ultimately, while
the exponent $E(N,3)$ of Lemma~\ref{lem:capbd} is essentially
$2N\log_3 N$, the actual exponent for this $\phi$
should be something more like $(4/3)N\log_3 N$.

In the archimedean case, the Chebyshev polynomials
$\{\psi_j\}_{j\geq 1}$ provide
an even stronger example of this phenomenon.
More precisely, if $\Cv=\CC$
and $\phi(z)=\psi_2(z)=z^2-2$,
then $\calK_v$ is simply the interval
$[-2,2]$ in the real line.  For $j\geq 1$,
the $j^{\text{th}}$
Chebyshev polynomial $\psi_j$ has $|\psi_j|\leq 2$ on $\calK_v$,
as compared with the proof's bound of $2^{c_0 + c_1 + \cdots}$
for $|f_j|$.

In general, however, knowing nothing about the polynomial
other than its degree and the radius $r_v$, we cannot substantially
improve on Lemma~\ref{lem:capbd}.
\end{remark}

\section{A Partition of the Filled Julia Set: Non-archimedean Case}
\label{sect:cap2}

The key to the Main Theorem,
as described in the introduction, is to divide the filled
Julia set at a particular bad prime into two smaller pieces
$X_1$ and $X_2$.  As a result, the
product $\prod_{i\neq j} |x_i-x_i|_v$,
when restricted to $\{x_i\}\subseteq X_k$ (for fixed
$k=1,2$), will be substantially smaller than the bound
of Lemma~\ref{lem:capbd}.  We begin with non-archimedean primes.

\begin{lemma}
\label{lem:napart}
Let $\Cv$ be a complete, algebraically closed field with
non-archimedean absolute value $|\cdot|_v$.
Let $\phi\in\Cv[z]$ be a
polynomial of degree $d\geq 2$ with lead coefficient
$a_d\in\Cv$.  Denote by $\calK_v$
the filled Julia set of $\phi$ in $\Cv$, and let $r'_v$
be the radius of the smallest disk $U_0$ that contains $\calK_v$.
Set $r_v = |a_d|_v^{1/(d-1)} r'_v$, and suppose
that $r_v>1$.

Then there are disjoint sets $X_1,X_2\subseteq\calK_v$
and positive integers $m_1,m_2$
with the properties
that $X_1\cup X_2=\calK_v$,
that $m_1+m_2=d$,
that for $k=1,2$,
$\phi:X_k\twoheadrightarrow \calK_v$ is $m_k$-to-$1$,
and that for $k=1,2$,
for any integer $N\geq 2$, and for any set $\{x_1,\ldots,x_N\}
\subseteq X_k$ of $N$ points in $X_k$,
$$\prod_{i\neq j} |x_i - x_j|_v \leq 
|a_d|_v^{-N(N-1)/(d-1)} r_v^{E(N,m_k,d)},$$
where $E(N,m_k,d)$ is as in Definition~\ref{def:E}
\end{lemma}

\begin{proof}
As in the proof of Lemma~\ref{lem:capbd}, we may assume
that $\phi$ is monic.

By Lemma~\ref{lem:goodrad}, $U_0$ is a closed disk
of radius $r_v\in |\Cv^{\times}|_v$.
We may write $U_0=\Dbar(a,r_v)$ for some point $a\in\calK_v$,
since $\calK_v$ is nonempty, and since any point of a non-archimedean
disk is a center.  Pick $b\in\phi^{-1}(a)$.  Note that
$b\in\calK_v\subseteq U_0$.

Write $U_1=\phi^{-1}(U_0)$.  By Lemma~\ref{lem:manydisk},
$U_1 = D_1 \cup\cdots \cup D_{\ell}$ for some disjoint
closed disks $\{D_i\}$, with $2\leq \ell \leq d$.
Moreover, $\phi: D_i\twoheadrightarrow U_0$ maps $d_i$-to-one
for some positive integers $\{d_i\}$ with $d_1+\cdots+d_{\ell}=d$.
Define 
$$W_1=\{x\in U_1 : |x-b|_v < r_v\},
\qquad\text{and}\qquad
W_2=U_1\setminus W_1,
$$
so that $W_1\cap W_2=\emptyset$
and $W_1\cup W_2=U_1$.  If $W_2=\emptyset$, then
$\calK_v\subseteq D(b,r_v)\subsetneq U_0$,
contradicting the minimality of $U_0$.  (The second inclusion
is strict because $r_v\in |\Cv^{\times}|_v$.)  Thus,
since $b\in W_1$, both $W_1$ and $W_2$ are nonempty.

Furthermore, $W_1$ and $W_2$ are both
finite unions of disks $D_i$ above.  Hence,
there are integers $m_1,m_2\geq 1$ so that each $W_k$ maps
$m_k$-to-one onto $U_0$, with $m_1+m_2=d$.  Let $X_k=W_k\cap\calK_v$
for $k=1,2$.  Since $\phi^{-1}(\calK_v)=\calK_v$,
$\phi$ must map $X_k$ $m_k$-to-one onto $\calK_v$.

For any integer $i\geq 1$,
observe that the polynomial $\phi^i(z) - a$ is monic of degree
$d^n$.  Moreover, since the equation $\phi^{i-1}(z)=a$ has exactly
$d^{i-1}$ roots (counting multiplicity), all of which lie in $U_0$,
it follows that $\phi^{i}(z)=a$ has $m_1 d^{i-1}$ roots in
$W_1$ and $m_2 d^{i-1}$ roots in $W_2$, counting multiplicity.
Thus, we may write
$$\phi^i(z) - a = g_i(z) h_i(z)$$
where $g_i$ is monic of degree $m_1 d^{i-1}$ with all its roots in $W_1$,
and 
$h_i$ is monic of degree $m_2 d^{i-1}$ with all its roots in $W_2$.
%In addition, define $g_0(z)=z-b_1$ and $h_0(z)=z-b_2$, where $b_1\in W_1$
%and $b_2\in W_2$.
In addition, define $g_0(z)=h_0(z)=z-a$.

We will now use the polynomials $g_i$ to compute the bounds given
in the Lemma for $X_1$; the proof for $X_2$ is similar, using $h_i$.
To simplify notation, write $X=X_1$ and $m=m_1$.

For any integer $j\geq 0$, write $j=c_0 + mk$, and write
$k$ in base-$d$ notation, so that
$$j= c_0 + m(c_1 + c_2 d + c_3 d^2 + \cdots + c_M d^{M-1}),$$
with $c_0\in\{0,1,\ldots, m-1\}$,
and with $c_i\in\{0,1,\ldots, d-1\}$ for $i\geq 1$.
Define
$$f_j(z) = \prod_{i=0}^M [g_i(z)]^{c_i}.$$
Clearly, $f_j$ is monic of degree $j$.
Meanwhile, for $x\in X$ and $i\geq 1$,
observe that $\phi^i(x)\in\calK_v$, and therefore
$|\phi^i(x)-a|\leq r_v$.  On the other hand, all roots of $h_i$
lie in $W_2$, which is distance $r_v$ from $x$; therefore,
$|h_i(x)| = r_v^{(d-m)d^{i-1}}$.  It follows that
$$|g_i(x)|_v \leq r_v^{1-(d-m) d^{i-1}}$$
for all $i\geq 1$.  In addition, since $X\subseteq U_0$,
we have $|g_0(x)| \leq r_v$.  Thus,
$$|f_j(x)|_v \leq r_v^{c_0} \prod_{i=1}^M r_v^{c_i(1-(d-m)d^{i-1})}
= r_v^e,$$
where
$e=e(j,m,d)$
%$$e = c_0 + (c_1 +\cdots + c_M) -(d-m)(c_1 + c_2 d + \cdots c_M d^{M-1})
%= c_0 + e(k,d) -(d-m)k = e(j,m,d),$$
in the notation of Definition~\ref{def:E}.

By the same van der Monde determinant argument as in
the proof of Lemma~\ref{lem:capbd}, it follows that
if $N\geq 2$ and $x_1,\ldots, x_N\in X$, then
\[
\prod_{i\neq j} |x_i - x_j|_v \leq r_v^{E(N,m,d)} .
\qedhere
\]
\end{proof}

\begin{remark}
\label{rem:sharp2}
In some cases,
$\calK_v$
splits into more than two pieces, each
much smaller than the $X_1,X_2$ of
Lemma~\ref{lem:napart}.  For example, the filled Julia set
of the map $\phi(z)=z^d -c^{d-1}z$ of Example~\ref{ex:cantor}
breaks naturally into $d$ pieces.
Adapting the method of the Lemma
for each piece, we could
ultimately replace the coefficient $d^2 -2d + 2$
in Theorem~\ref{thm:global} by $d$.

However, as previously noted, most polynomials are not so simple.
Indeed, the filled Julia set
of $\phi(z) = z^d - a z^{d-1}$ from Example~\ref{ex:hsia}
splits into only two pieces.  (Of course, if we take
a higher preimage $U_n$ in that example, we get
more than two pieces; but because
of the large radii, there appears to
be no improvement gained by using $n>1$.)
Even an application of 
the arguments of Remark~\ref{rem:sharp1}
would result in
only a slight decrease in the
coefficient of $N\log_d N$ in the exponent
(cf.~Lemma~\ref{lem:EFbd}.b).
Unfortunately, a real improvement would require
an increase in the size of the
(negative) coefficient of $N^2$, not the $N\log_d N$ term.
\end{remark}

\section{A Partition of the Filled Julia Set: Archimedean Case}
\label{sect:cap3}

The final tool needed for Theorem~\ref{thm:global} is
an archimedean analogue of Lemma~\ref{lem:napart}.  Roughly
the same argument works, but only if the diameter of
the filled Julia set $\calK$ is large enough.  This phenomenon
is familiar to complex dynamicists.
For example, given $\phi(z) = z^2 + c \in\CC[z]$,
if the diameter of $\calK$ is
small, then $c$ lies in the Mandelbrot
set, in which case $\calK$ is connected.  However,
once the diameter is large enough, $c$ leaves the Mandelbrot
set and $\calK$ becomes disconnected.  In fact,
as the diameter grows, the various
pieces of $\calK$ shrink.

We begin with the following preliminary result.

\begin{lemma}
\label{lem:preapart}
Let $\phi\in\CC[z]$ be a
polynomial of degree $d\geq 2$ with lead coefficient
$a_d\in\CC$.  Denote by $\calK$
the filled Julia set of $\phi$ in $\CC$,
and let $U_0=\Dbar(a,r')$
be the smallest disk that contains $\calK$.
Set $r = |a_d|^{1/(d-1)} r'$, and suppose that
$$r >
\begin{cases}
3, & \text{if } d=2, \text{ or}
\\
2 + \sqrt{3}, & \text{if } d\geq 3.
\end{cases}
$$
Then $\calK$ is contained in the union of $d$
open disks of radius $|a_d|^{-1/(d-1)}$.
\end{lemma}

\begin{proof}
As in the proof of Lemma~\ref{lem:capbd}, we may assume
that $\phi$ is monic.
Denote by $b_1,\ldots, b_d$ the (possibly repeated) roots
of $\phi(z)=a$, and let $\Dbar(c,s)$ be the smallest disk
containing $b_1,\ldots,b_d$.
(Here, we break our convention and allow $s=0$ if
$b_1=\cdots=b_d$.)
Because $\calK$ is not contained in $D(c,r)$, there must
be some $y_0\in\calK$ such that $|y_0-c|\geq r$.

Let $Y = \Dbar(c,s) \cap D(y_0,|y_0-c|)$.
We claim that $Y$ is contained
in a disk of radius strictly less than $s$
(or that $Y$ is empty, if $s=0$).
Indeed,
if $|y_0-c|<s$, then $Y\subseteq D(y_0,|y_0-c|)$ trivially.
Otherwise, $|y_0-c|\geq s$, and since the center $c$ of the
first disk lies on the boundary of the second,
the intersection $Y$ is contained in a strictly smaller disk.
(For example, center the new disk at the midpoint of the
two intersection points of the two boundary circles.)

By the minimality of $s$, then, not all of
$b_1,\ldots,b_d$ can be in $Y$.
Thus, there is some $1\leq i\leq d$
such that $|y_0-b_i|\geq |y_0-c|\geq r$.
Without loss, assume that $|y_0-b_1|\geq r$.

For all $i\geq 2$, we have $|y_0-b_i|\geq r-s$, because $b_i\in\Dbar(c,s)$.
Since $y_0\in\calK$, we have $\phi(y_0)\in\calK$, and therefore
$|\phi(y_0)-a|\leq r$.  If $r-s\geq 0$, then, we have
$$r\geq |\phi(y_0)-a| = \prod_{i=1}^d |y_0-b_i|
= |y_0-b_1| \cdot \prod_{i=2}^d |y_0-b_i|
\geq r\cdot (r-s)^{d-1},$$
from which we obtain $r-s\leq 1$.  Regardless of the sign of $r-s$, then,
we have $s\geq r-1$.

Re-index so that $b_1$ and $b_d$ are distance $\max\{|b_i-b_j|\}$
apart, and so that for all $i=1,\ldots, d-1$, we
have $|b_{i+1}-b_1| \geq |b_i - b_1|$.  Thus,
the $\{b_i\}$ are ordered by their distance from $b_1$.
Moreover,
$|b_d - b_1|\geq\sqrt{3} s \geq \sqrt{3} (r-1)$; see, for
example, \cite{YB}, Exercise~6-1.

If $d=2$, we can improve this lower bound.  In that case,
the smallest disk containing $b_1$ and $b_2$
is the closed disk centered
at $(b_1+b_2)/2$ of radius $|b_1-b_2|/2$.  That is,
$s=|b_1-b_2|/2$.  It follows that
$|b_1-b_2| =2s\geq 2(r-1)$.

For all degrees $d\geq 2$, we have $r>3$, so
that $s>2$, and therefore the two disks
$\Dbar(b_1,1)$ and $\Dbar(b_d,1)$ are disjoint.
Moreover,
as $y$ ranges through $\CC\setminus[D(b_1,1)\cup D(b_d,1)]$,
the minimum value of $|y-b_1|\cdot|y-b_d|$ is $|b_1-b_d|-1$,
attained at only two points, namely the point on the boundary
of each disk closest to the other disk.

Let $U_1 = \phi^{-1}(U_0)$.
Since $\calK=\phi^{-1}(\calK)\subseteq U_1$,
it suffices to show that
\begin{equation}
\label{eq:Uinc}
U_1 \subseteq \bigcup_{i=1}^d D(b_i,1).
\end{equation}
If not, then
there is some $y\in U_1 \setminus \bigcup D(b_i,1)$.
If $d\geq 3$, then by the above computations, we have
$|y-b_1|\cdot |y-b_d| \geq \sqrt{3}(r-1) -1$.
Since $\phi(y)\in U_0$, we obtain
$$r\geq |\phi(y)-a| = \prod_{i=1}^d |y-b_i|
\geq (\sqrt{3}(r-1) -1) \prod_{i=2}^{d-1} |y-b_i|
\geq (\sqrt{3}(r-1) -1),$$
contradicting the hypothesis that $r>2+\sqrt{3}$.  Similarly,
if $d=2$, then
$$r\geq |\phi(y)-a| = |y-b_1|\cdot |y-b_2| > 2(r-1) - 1
 = 2r -3,$$
contradicting the hypothesis that $r> 3$,
and proving the Lemma.
\end{proof}

\begin{remark}
\label{rem:Kcomp}
Because $\calK_v$ is compact for archimedean $v$,
the conclusion of Lemma~\ref{lem:preapart} implies
that $\calK$ is in fact contained in $d$ closed disks
of radius strictly less than $|a_d|_v^{-1/(d-1)}$.  This fact
will be useful in Cases~2 and~3 of the proof of
Theorem~\ref{thm:global}.
\end{remark}

We are now prepared to present our archimedean version
of Lemma~\ref{lem:napart}.

\begin{lemma}
\label{lem:apart}
Let $\phi\in\CC[z]$ be a
polynomial of degree $d\geq 2$ with lead coefficient
$a_d\in\CC$.  Denote by $\calK$
the filled Julia set of $\phi$ in $\CC$, and let $r'$
be the radius of the smallest disk $U_0$ that contains $\calK$.
Set $r = |a_d|^{1/(d-1)} r'$ and
$$C_d = d^{-(d-2)/(d-1)}
\leq \min\left\{1, \frac{1.2}{d-1} \right\}.$$
Suppose that
$$r \geq 
\begin{cases}
4
& \text{if } d=2
\\
\dfrac{\sqrt{3} + 2(d-1)}{\sqrt{3} - (d-1)C_d},
& \text{if } d\geq 3.
\end{cases}
$$

Then there are disjoint sets $X_1,X_2\subseteq\calK$
and positive integers $m_1,m_2$
with the properties
that $X_1\cup X_2=\calK$,
that $m_1+m_2=d$,
that for $k=1,2$,
$\phi:X_k\twoheadrightarrow \calK_v$ is $m_k$-to-$1$,
and that for $k=1,2$,
for any integer $N\geq 2$, and for any set $\{x_1,\ldots,x_N\}
\subseteq X_k$ of $N$ points in $X_k$,
$$\prod_{i\neq j} |x_i - x_j| \leq
N^N |a_d|^{-N(N-1)/(d-1)}
C_d^{-F(N,m_k,d)}\left(C_d r\right)^{E(N,m_k,d)},$$
where $E(N,m_k,d)$ and $F(N,m_k,d)$ are as in Definition~\ref{def:E}.
\end{lemma}

\begin{proof}
As in the proof of Lemma~\ref{lem:capbd}, we may assume
that $\phi$ is monic.
It is easy to check that $C_d \leq \min\{ 1, 1.2/(d-1)\}$
(the closest approach for $d\geq 3$ occurs at $d=5$),
and that
the lower bound $(\sqrt{3} + 2(d-1))/(\sqrt{3} - (d-1)C_d)$
(respectively, $4$)
for $r$ is greater than $2 +\sqrt{3}$ (respectively, $3$),
so that we may invoke Lemma~\ref{lem:preapart}.

Write $U_0=\Dbar(a,r)$, and
define and order
$b_1,\ldots,b_d$ as in the proof of Lemma~\ref{lem:preapart},
so that $|b_1-b_d| \geq \sqrt{3} (r-1)$ (or
$|b_1-b_d| \geq 2(r-1)$, if $d=2$).

If $d\geq 3$, observe that for some $m=1,\ldots,d-1$,
we have
$$|b_{m+1} - b_1| \geq |b_m - b_1| + 2 +  C_d r.$$
For if not, then
$$\sqrt{3}(r-1) \leq |b_d - b_1| < (d-1)\left[ 2 + C_d r \right]
= 2(d-1) + (d-1) C_d r,$$
so that $[\sqrt{3} - (d-1) C_d ] r < \sqrt{3} + 2(d-1)$, contradicting
the hypotheses.

If $d=2$, we have $|b_2-b_1|\geq 2r -2 \geq 2 + r$, since $r\geq 4$.
Let $m=1$ in this case.

Let $U_1=\phi^{-1}(U_0)$, and set
$W_1 = \Dbar(b_1,|b_m - b_1| + 1) \cap U_1$
and $W_2 = U_1\setminus W_1$.  Observe
that $\dist(W_1,W_2)\geq C_d r$.
Indeed, if $y_1\in W_1$ and $y_2\in W_2$,
then $y_1\in \Dbar(b_i,1)$ and $y_2\in \Dbar(b_j,1)$ for
some $1\leq i \leq m$ and some $m+1\leq j \leq d$; therefore
$$|y_2 - y_1| \geq |b_j- b_1| - |b_i - b_1| - 2
\geq |b_{m+1} - b_1| - |b_m - b_1| - 2 \geq C_d r.$$
Since $W_1$ contains $m$ preimages of $a$ and $W_2$ contains
the other $d-m$, it follows that $\phi$ maps $W_1$ $m$-to-$1$
onto the connected set $U_0$, and it maps $W_2$ $(d-m)$-to-$1$
onto $U_0$.

Let $X_1=W_1\cap\calK$, $X_2=W_2\cap\calK$, $m_1=m$,
and $m_2=d-m$.  By the previous paragraph, $X_1$ and $X_2$
satisfy all of the mapping properties claimed in the Lemma.
For any integer $i\geq 0$, define $g_i(z)$ and $h_i(z)$
as in the proof of Lemma~\ref{lem:napart}.  That is,
for $i\geq 1$, write
$$\phi^i(z) - a = g_i(z) h_i(z),$$
where $g_i$ is a monic polynomial of degree $m_1d^{i-1}$
with all of its roots in $W_1$, and $h_i$ is a monic
polynomial of degree $m_2 d^{i-1}$ with all of its roots
in $W_2$.  For $i=0$, define $g_0(z)=h_0(z) = z-a$.
We will now compute the bounds given in the Lemma for $X_1$;
the proof for $X_2$ is similar.  Write $X=X_1$ and $m=m_1$.

As in the proof of Lemma~\ref{lem:napart}, we may write
any integer $j\geq 0$ as
$$j= c_0 + m(c_1 + c_2 d + c_3 d^2 + \cdots + c_M d^{M-1}),$$
with $c_0\in\{0,1,\ldots, m-1\}$,
and with $c_i\in\{0,1,\ldots, d-1\}$ for $i\geq 1$.  Similarly,
define
$$f_j(z) = \prod_{i=0}^M [g_i(z)]^{c_i},$$
which is clearly monic of degree $j$.  As before, for any $x\in X$, we have
$|\phi^i(x)-a|\leq r$.  Similarly, the roots of $h_i$,
which all lie in $W_2$ (for $i\geq 1$),
are distance at least $C_d r$ from $x$
(which is worse than the Lemma~\ref{lem:napart} bound of $r$).  Thus,
$|h_i(x)| \geq (C_d r)^{(d-m)d^{i-1}}$, and hence
$$|g_i(x)| \leq r (C_d r)^{-(d-m) d^{i-1}}
= C_d^{-1} (C_d r)^{1-(d-m) d^{i-1}}$$
for all $i\geq 1$.  Moreover, since $x\in U_0$, we have $|g_0(x)|\leq r$.
We obtain
$$|f_j(x)| \leq r^{c_0} \prod_{i=1}^M C_d^{-c_i} (C_d r)^{c_i(1-(d-m)d^{i-1})}
= C_d^{-(c_0 + c_1 + \cdots + c_M)} (C_d r)^e,$$
where
$e=e(j,m,d)$
%$$e = c_0 + (c_1 +\cdots + c_M) -(d-m)(c_1 + c_2 d + \cdots c_M d^{M-1})
%= c_0 + e(k,d) -(d-m)k = e(j,m,d),$$
in the notation of Definition~\ref{def:E}.  The Lemma then
follows by the van der Monde determinant argument of
the proof of Lemma~\ref{lem:capbd}.
\end{proof}

\begin{remark}
\label{rem:cutoff}
Later, in the proof of Theorem~\ref{thm:global}, we will
consider the quantity $C_d r$, rather than the radius $r$,
at the archimedean primes.  It is easy
to prove that the lower bound for $r$
given in Lemma~\ref{lem:apart} is guaranteed to hold
provided $C_d r \geq 4 + \sqrt{3}$.  (In fact, $4+\sqrt{3}$ is
the value of $C_d(\sqrt{3} + 2(d-1))/(\sqrt{3} - (d-1)C_d)$
at $d=3$.)  For $d=2$, we also note the more obvious facts
that $C_2=1$ and that
the corresponding sufficient lower bound for $C_2 r$ is $4$.
\end{remark}

\begin{remark}
\label{rem:archbad2}
The bounds in Lemma~\ref{lem:preapart} and
Lemma~\ref{lem:apart} are not sharp.
Besides the fact that 
most of the comments
from Remarks~\ref{rem:sharp1} and~\ref{rem:sharp2}
apply here, our geometric
arguments could also be improved.
For example, in the proof of Lemma~\ref{lem:preapart}, if
we considered $\Dbar(c,s)\cap D(y_0,t)$ instead
of $\Dbar(c,s) \cap D(y_0,|y_0-c|)$, where
$t=\sqrt{|y_0-c|^2 + s^2}$, we
could show that some $b_i$ satisfies
$|y_0-b_i|\geq t$.  Related arguments could show that two or
more points $b_i,b_j$ must make the product $|y_0-b_i|\cdot |y_0-b_j|$
larger than we proved.  Similarly, it should be possible to
increase the
$\sqrt{3}$ factor to something closer to $2$
by considering the geometric arrangement of the $\{b_i\}$
more delicately.
\end{remark}

\section{The Global Bound}
\label{sect:global}

At last, we are prepared to state and prove a precise version of
the Main Theorem.

\begin{thm}
\label{thm:global}
Let $K$ be a global field, and let $\phi\in K[z]$ be a polynomial
of degree $d\geq 2$.  
Let $s_{\infty}\geq 0$ be the number of archimedean primes of $K$,
and let $s\geq s_{\infty}$
be the number of bad {\em (}i.e., not potentially good{\em )}
primes of $\phi$ in $M_K$, including
all archimedean primes.

If $K$ is a function field, let $q$ be the size of the smallest
residue field of a prime $v\in M_K$.
If $K$ is a number field, let $D=[K\col\QQ \hspace{0.75pt}]$,
and let
$$
\sigma = \begin{cases}
7 & \text{if } d=2, 
\\
\dfrac{2\cdot 33^{(d-1)(d-2)}}{(d-1)(d-2)} & \text{if } d\geq 3.
\end{cases}
$$ 
Set
$$t =
\begin{cases}
s - s_{\infty}
& \text{if $K$ is a number field and } s \leq \sigma D ,
\\
s + \dfrac{D \log d}{2\log 2}
& \text{if $K$ is a number field and } s > \sigma D,
\\
s & \text{if $K$ is a function field,}
\end{cases}
$$
and
$$
\beta = \begin{cases}
9 & \text{if $K$ is a number field, } s \leq \sigma D,
\text{ and }  d=2, 
\\
\max\{11,2d\} & \text{if $K$ is a number field, } s \leq \sigma D,
\text{ and }  d\geq 3,
\\
1 & \text{otherwise.}
\end{cases}
$$ 
Then $\phi$ has no more than $M+1$ $K$-rational preperiodic points
in $\PK$, where
$$M =
\begin{cases}
q & \text{if $K$ is a function field and } s=0,
\\
\beta^D & \text{if $K$ is a number field and } s=s_{\infty},
\\
\beta^D(d^2 - 2d + 2)(t\log_d t + 3t)
& \text{if } 0<t<d,
\\
\beta^D(d^2 - 2d + 2)(t\log_d t + t\log_d\log_d t + 3t)
& \text{otherwise.}
\end{cases}
$$
%In particular, for fixed degree $d\geq 2$ and
%fixed $D$ (or fixed $q$, as appropriate),
%the number of $K$-rational preperiodic points of $\phi$
%is $O(s\log s)$, where $s$ is the number of bad primes.
\end{thm}

\begin{proof}
For each prime $v\in M_K$, let $n_v\geq 1$ be the
exponent so that the product formula
\eqref{eq:pfmla}
holds for all $x\in K^{\times}$.
Let $S$ be the (finite) set of primes of $K$ of bad reduction
of $\phi$,
including all the archimedean primes; that is, $\#S = s$.
Let $a_d\in K$ be the lead coefficient of $\phi$.
For each prime $v\in M_K$, let $\calK_v\subseteq\Cv$ denote
the filled Julia set of $\phi$ in $\Cv$,
let $r'_v$ be the radius of
the smallest disk in $\Cv$ containing $\calK_v$,
and let $r_v = |a_d|_v^{1/(d-1)} r'_v$.

For each non-archimedean prime $v$, let $R_v = r_v^{n_v}$.
For each archimedean prime $v$, let $R_v = (C_d r_v)^{n_v}$,
where $C_d = d^{-(d-2)/(d-1)}\leq 1$, as in the statement
of Lemma~\ref{lem:apart}.
We consider four cases, some of which overlap with others.

{\bf Case 0.}  The simplest case is that $K$ is a function
field and $S=\emptyset$; that is, there are no archimedean
primes, and all primes have potentially good reduction.
Let $w\in M_K$ be a prime whose residue field has only
$q$ elements, and
suppose that there are $q+1$ distinct $K$-rational preperiodic
points $\{x_1,\ldots,x_{q+1}\}$ besides the point at $\infty$.

By Lemma~\ref{lem:goodrad}.c, we have
$|x_i - x_j|_v \leq |a_d|_v^{-1/(d-1)}$ for every
$v\in M_K$ and every $i,j\in\{1,\ldots,n\}$.
Moreover, by the pigeonhole principle, there must
be some distinct $i,j\in\{1,\ldots,n\}$ such that
$|x_i - x_j|_w < |a_d|_w^{-1/(d-1)}$.
Hence,
$$
1 = \prod_{v\in M_K} \left|x_i - x_j\right|_v^{n_v}
< \prod_{v\in M_K} \left[ |a_d|_v^{-1/(d-1)}
  \right]^{n_v}
= 1,
$$
which is a contradiction.
Thus, there are at most $q$ finite $K$-rational 
preperiodic points.

{\bf Case 1.}
Choose $w\in M_K$ such that
$R_w \geq R_v$ for all $v\in M_K$.
(Such $w$ exists because $R_v=1$ for all but finitely many
$v\in M_K$.)
In this main case, we suppose that:
\begin{itemize}
\item $R_w > 1$.
\item If $K$ is a number field, then $R_w\geq 4$ and $s-s_{\infty}\geq 1$.
\item If $w$ is archimedean, then the lower bounds
  of Lemma~\ref{lem:apart} hold for $r_w$.
\end{itemize}
In particular, we may choose integers $m_1, m_2$
and sets $X_1,X_2\subseteq \calK_w$ for $\phi$ according
to Lemma~\ref{lem:napart} (if $w$ is non-archimedean)
or Lemma~\ref{lem:apart} (if $w$ is archimedean).

For each index $k=1,2$, set
$$A_k = \frac{d-m_k}{m_k(d-1)},
\qquad
B_k = 1-\log_d m_k,
\qquad\text{and}\qquad
N_k = M(A_k,B_k,t),$$
where $M(\cdot,\cdot,\cdot)$ is as in
Lemma~\ref{lem:slogs}, and where
$t$ is as in the statement of the Theorem.
We claim that
there are fewer than $N_k$ $K$-rational preperiodic
points in $X_k$.

To prove the claim, fix $k=1,2$, and let $m=m_k$, $A=A_k$,
$B=B_k$, and $N=N_k$.
Suppose there are $N$ distinct $K$-rational preperiodic points
$x_1,\ldots, x_N$ in $X_k$.
Then by the product formula
applied to both $\prod_{i\neq j} (x_i-x_j)$ and $a_d$,
\begin{align}
\label{eq:mainprod}
1 &= \prod_{v\in M_K} \left|\prod_{i\neq j} (x_i - x_j)\right|_v^{n_v}
= \prod_{v\in M_K} \left[ |a_d|_v^{N(N-1)/(d-1)}
  \prod_{i\neq j} |x_i - x_j|_v\right]^{n_v}
\notag
\\
& \leq  \prod_{v\in S} \left[ |a_d|_v^{N(N-1)/(d-1)}
  \prod_{i\neq j} |x_i - x_j|_v \right]^{n_v},
\end{align}
where the inequality is because $|x-y|_v\leq |a_d|_v^{-1/(d-1)}$ for all
$v\in M_K\setminus S$ and $x,y\in \calK_v$, by Lemma~\ref{lem:goodrad}.c,
and because $x_1,\ldots,x_N\in\calK_v$ for every $v\in M_K$.

If $w$ is non-archimedean, then by Lemma~\ref{lem:capbd}
and Lemma~\ref{lem:napart}, \eqref{eq:mainprod} becomes
$$ 1 \leq N^{DN} r_w^{n_w E(N,m,d)}
\prod_{v\in S\setminus\{w\} } r_v^{n_v E(N,d)}
= N^{DN} C_d^{-D E(N,d)} R_w^{E(N,m,d)}
\prod_{v\in S\setminus\{w\} } R_v^{E(N,d)},$$
where we set $D=0$ if $K$ is a function field.
(The appearance of $D$ in the exponent comes
from equation~\eqref{eq:sumnv}.)
Because $R_w\geq R_v$ and $E(N,d)\geq 0$, we can
replace each $R_v$ by $R_w$; and
because $R_w, C_d^{-1}\geq 1$, we can apply
Lemma~\ref{lem:EFbd}.a--b to obtain
\begin{equation}
\label{eq:midbound}
1 \leq N^{DN} C_d^{-(d-1)D N\log_d N}
R_w^{(d-1)N[ s\log_d N - AN + B ]}.
\end{equation}
Similarly, if $w$ is archimedean, then by Lemma~\ref{lem:capbd}
and Lemma~\ref{lem:apart}, \eqref{eq:mainprod} becomes
\begin{align*}
1 & \leq N^{DN} C_d^{-n_w F(N,M,d)} (C_d r_w)^{n_w E(N,m,d)}
\prod_{v\in S\setminus\{w\} } r_v^{n_v E(N,d)}
\\
& = N^{DN} C_d^{-(D-n_w) E(N,d) - n_w F(N,m,d)} R_w^{E(N,m,d)}
\prod_{v\in S\setminus\{w\} } R_v^{E(N,d)}.
\end{align*}
Replacing each $R_v$ by $R_w$ as before, and applying
Lemma~\ref{lem:EFbd}a--c,
we obtain exactly inequality~\eqref{eq:midbound} once more.

Meanwhile, we compute
\begin{equation}
\label{eq:facbound}
N^{DN} C_d^{-(d-1)D N\log_d N}
= d^{DN \log_d N} d^{(d-2)(DN\log_d N)} = d^{(d-1)DN\log_d N}.
\end{equation}
%If $K$ is a number field, then our assumed lower bounds
%for $r_w$ in this case imply that $R_w\geq 4$.
%(In fact, for $d\geq 3$, the bound $C_d\geq 1/(d-1)$
%shows that $R_w\geq 2(\sqrt{3}+1)$.)  Thus,
If $K$ is a number field, our assumption that $R_w\geq 4$
means that
$d\leq R_w^{1/\log_d 4}$.
Combining \eqref{eq:midbound} and \eqref{eq:facbound}, then, we
obtain
\begin{equation}
\label{eq:mainleq}
1\leq R_w^{(d-1)N [ t \log_d N - AN + B]},
\end{equation}
where $t=s + D \log d/(2\log 2)$, as in the statement of
the Theorem.  The same inequality follows for function fields
with $t=s$, since $D=0$ in that case.
By our definitions of $A$, $B$, and $t$,
the hypotheses of Lemma~\ref{lem:slogs} hold.
Thus, by that Lemma and our choice
of $N$, we have $t \log_d N - AN + B < 0$, so that $1<1$,
which is a contradiction, proving the claim that
there are fewer than $N_k$ $K$-rational preperiodic
points in $X_k$.  (However, since $N_k$ need not be
an integer, we cannot claim that there are at most
$N_k - 1$ such points.)

The total number of finite $K$-rational preperiodic points is
the number in $X_1$ plus the number in $X_2$.  That is,
there are fewer than $N_1 + N_2$ such points.
That upper bound is
\begin{equation}
\label{eq:Nsum}
N_1 + N_2 = M(A_1,B_1,t) + M(A_2,B_2,t).
\end{equation}
From the definition of $M(A,B,t)$ in
Lemma~\ref{lem:slogs}, it is each to check that, as $m_1$ varies
from $1$ to $d-1$, the largest value of $N_1 + N_2$
in equation~\eqref{eq:Nsum} is attained at
$m_1=1$ and $m_2=d-1$ (or vice versa).  In that case, the bound is
$$N_1 + N_2 = 
(d^2 - 2d + 2)
\left(t \log_d t + t \log_d(\max\{1,\log_d t\}) + 3t\right) .
$$
Adding $1$ for the point at $\infty$,
we obtain the bound stated in the Theorem, with $\beta=1$.

{\bf Case 2.}  Next, suppose that $K$ is a number field and $d=2$.
Write $S_{\infty}$ for the set of archimedean primes of $M_K$,
and let $s_{\infty} = \# S_{\infty}$.
We will remove the archimedean primes from
the picture by covering the filled Julia set at each such prime
$v\in S_{\infty}$
by at most $9^{n_v}$ disks of diameter less than $|a_d|_v^{-1}$.
To simplify notation, let
$\calK'_v = a_d \calK_v$; we wish
to cover $\calK'_v$ by disks of diameter less than $1$.

For any real prime $v\in S_{\infty}$, the set
$\calK_v'$ is contained either in a single interval of length $6$
or in two intervals of length less than $2$,
by Lemma~\ref{lem:preapart} and Remark~\ref{rem:Kcomp}.
(In fact, the bound of $6$ could be reduced to $4$,
but we will not need that stronger bound here.)
In particular, $\calK'_v$ is contained in a union of
seven or fewer intervals of length strictly less than $1$.

For a complex prime $v\in S_{\infty}$, the same Lemma implies that
$\calK'_v$ is contained either in a single disk of radius
$3$ or in two disks of radius less than $1$.
Each disk of radius $1$
can easily be covered by nine disks of diameter slightly less
than $1$.  Similarly,
the disk of radius $3$
can be covered by a square of side length $6$.
That square can then be divided into $81$ squares
of side length $2/3$, each
of which fits inside a disk of diameter less than $1$.

Scaling back by $|a_d|_v^{-1}$,
then, we have at each archimedean prime $v\in S_{\infty}$
at most $9^{n_v}$ disks of diameter less
than $|a_d|_v^{-1}$ which together cover $\calK_v$,
as promised.
In total, then, we have at most $9^D$ choices of one disk for each
archimedean prime.

For any such choice $\DD=\{D_v:v\in S_{\infty}\}$ of one disk
of diameter less than
$|a_d|_v^{-1}$ for each archimedean prime $v$, let $\calP_{\DD}$
denote the set of $K$-rational preperiodic points $x$ for which
$x\in D_v$ for every $v\in S_{\infty}$.  We will bound the size
of $\calP_{\DD}$.

If $S=S_{\infty}$, then each set $\calP_{\DD}$ can contain
at most one point.  Indeed, if there were distinct points
$x,y\in\calP_{\DD}$, then
$$
1 = \prod_{v\in M_K} |x-y|_v^{n_v}
< \prod_{v\in M_K} \left[ |a_d|_v^{-1}
  \right]^{n_v}
= 1,
$$
by Lemma~\ref{lem:goodrad}.c, with the strict inequality
coming from the fact that the diameter at each archimedean
prime is strictly less than $|a_d|_v^{-1}$.
Since there are $9^D$ choices of $\DD$, there
are at most $9^D$ finite $K$-rational preperiodic points.

On the other hand, if $S\supsetneq S_{\infty}$, then
choose $w\in M_K\setminus S_{\infty}$ such that
$R_w \geq R_v$ for all $v\in M_K\setminus S_{\infty}$.
By Lemma~\ref{lem:goodrad}.c, $r_w>1$, so that we may
apply Lemma~\ref{lem:napart} at $w$.

Now fix $\DD$ and follow the argument of Case~1, but restricted to
$\{x_i\}\subseteq \calP_{\DD}$.  At each archimedean prime
$v\in S_{\infty}$ we have $|x_i-x_j|_v\leq |a_d|_v^{-1}$.  Therefore,
by Lemmas~\ref{lem:capbd}, \ref{lem:napart},
and~\ref{lem:EFbd}.a--b, inequality
\eqref{eq:mainprod} becomes
$$
1 \leq R_w^{(d-1)N[ t\log_d N - AN + B ]},
$$
where $t=s-s_{\infty}$.
Following the rest of
the argument of Case~1 (from inequality \eqref{eq:mainleq} on),
and multiplying by $9^D$ (the number of choices $\DD$),
we obtain the desired bounds.

{\bf Case 3.}
If $K$ is a number field and $d\geq 3$,
we proceed roughly as in Case~2.
Again, write $S_{\infty}$ for the set of archimedean primes of $M_K$,
let $s_{\infty} = \# S_{\infty}$, and let
$\calK'_v = \alpha \calK_v$,
where $\alpha^{d-1}=a_d$.
This time, we will cover $\calK'_v$
by at most $\beta^{n_v}$ disks of diameter less than $1$,
where $\beta = \max\{11, 2d\}$.

For a real prime $v\in S_{\infty}$, Lemma~\ref{lem:preapart}
and Remark~\ref{rem:Kcomp}
imply that $\calK'_v$
is contained either in a single interval of length
$4+2\sqrt{3}$
or in $d$ intervals of length less than $2$.
In particular, $\calK'_v$ is contained in a union of
$\max\{8,2d\}\leq \beta$ or fewer intervals of length
less than $1$.

For a complex prime $v\in S_{\infty}$, the same Lemma implies that
$\calK_v$ is contained either in a single disk of radius
$2+\sqrt{3}$
or in $d$ disks of radius less than $1$.  As before,
each disk of radius $1$
can be covered by nine disks of diameter less than $1$.
Similarly, the disk of radius $2+\sqrt{3}$
can be covered by a square of side length
$4+2\sqrt{3}$.
That square
can be divided into $121$
squares of side length $(4+2\sqrt{3})/11$,
each of which fits inside a disk of diameter less than $1$.
(In fact, using a hexagonal tiling, one could cover the big disk by
$84$ disks of diameter less than $1$,
but the messy proof
gives only a minor improvement over $121$.)
Thus, $\calK'_v$ can be covered by a union of
$\max\{121,9d\}\leq \beta^2$ disks of diameter less than $1$.

The rest of Case~3 then follows Case~2, with $\beta^D$ in place
of $9^D$.  This completes our analysis of the four cases.

{\bf Final step.}
If $K$ is a function field, we are done; indeed,
by Lemma~\ref{lem:goodrad}.c, Cases~0 and~1 cover all
the possibilities.

If $K$ is a number field, we will now show that
for $s>\sigma D$, we are automatically in Case~1.
Because $n_v\leq 2$ for an archimedean prime $v$,
and by Remark~\ref{rem:cutoff}, we need
only show there is some $w\in M_K$ such that $R_w \geq 4^2$ if $d=2$,
or such that $R_w \geq (4+\sqrt{3})^2$ if $d\geq 3$.

From basic algebraic number theory, there are
at most $D$ primes of $K$ above any given prime of $\QQ$.
Given an integer $m\geq 1$, let $p_m$ denote the $m^{\text{th}}$
prime in $\QQ$.  (That is, $p_1=2$, $p_2=3$, $p_3=5$, and so on.)
Thus, if $s-s_{\infty} >  D(m-1)$, there must be some
$w\in S\setminus S_{\infty}$ lying above a prime $p\geq p_m$ of $\QQ$.
Since $D\geq s_{\infty}$, we get such a $w$ provided
$s> mD$.

Meanwhile, by Lemma~\ref{lem:minrad}, given $w\in S\setminus S_{\infty}$
lying over a prime $p$ of $\QQ$, we have
$$R_w \geq  |\pi_w|_w^{-n_w} \geq p
\qquad \text{if } d=2,$$
or 
$$R_w^{(d-1)(d-2)} \geq |\pi_w |_w^{-n_w} \geq p
\qquad \text{if } d\geq 3,$$
where $\pi_w$ is a uniformizer at $w$.  (The Lemma applies
because if $\calK_w \cap K=\emptyset$, then there are no
finite $K$-rational preperiodic points at all, and the
conclusion of the Theorem is trivial.)
For $d=2$, then, the condition $R_w\geq 16$ is guaranteed
provided $s\geq 7D+1$, since $17$ is the seventh prime of $\QQ$.
Thus, $s > 7D = \sigma D$ suffices for $d=2$.

For $d\geq 3$,
the elementary estimate in
Theorem~4.7 of \cite{Apo} says that
$p_m > (1/6)m\log m$
for any integer $m\geq 1$.
It is easy to check that $m=\lfloor \sigma \rfloor$
satisfies $m\log m\geq 6(4+\sqrt{3})^{2(d-1)(d-2)}$,
where $\sigma= 2\cdot 33^{(d-1)(d-2)}/[(d-1)(d-2)]$
as in the statement of the Theorem.
(The $33$ appears because it is the smallest integer larger
than $(4+\sqrt{3})^2$.)
Thus, $s > \sigma D$ implies $R_w\geq (4+\sqrt{3})^2$,
once again forcing Case~1.
\end{proof}

\begin{remark}
If, for a given polynomial $\phi$, we know that we are in Case~1
(say, by inspection of the filled Julia set at one prime),
then we can set $\beta=1$ in the statement of the Theorem,
even if $s\leq\sigma D$.
In particular, for a fixed function $\phi$,
the conditions of Case~1 are preserved
if one passes to a finite extension of $K$.
Thus, one would
not have to worry about the growth of $s$ relative to $\sigma D$
as one traveled up a tower of number fields,
even though one cannot expect $s$ to increase as fast as $\sigma D$
in general.
\end{remark}

\begin{remark}
\label{rem:improve}
Our covering methods in Cases~2 and~3 are rather crude,
and it should be possible to cover the filled Julia sets $\calK_v$
in these cases far more efficiently.  For example,
our use of disks of diameter $1$ was rather simplistic.
Instead, one could cover $\calK_v$ by larger sets $Y$ for which
$\prod_{1\leq i,j\leq L} |y_i - y_j|_v \leq 1$ for some fixed
small integer $L$.  Such modifications could
lead to a substantial reduction in the coefficient $\beta^D$
in the final bound.  

Even without any extra work, the coefficient
can be improved in special cases.  For example,
if $K$ is a totally real number field, then the cutoff
$\sigma$, which determines when $\beta$ drops to $1$,
would be much smaller, since we would only need $R_w\geq 4$
(if $d=2$), or $R_w\geq 4 + \sqrt{3}$ (if $d\geq 3$)
rather than $4^2$ or $(4+\sqrt{3})^2$.
Moreover, if $K$ is totally real and $d=2$, then each archimedean
filled Julia set is contained in a union of four intervals
of length $1$.  (See, for example, Lemma~6.4 and Proposition~6.6
of \cite{CaGol}.)
Thus, the coefficient of $9^D$ that appears
in Theorem~\ref{thm:global} could be replaced by $4^D$,
with one exception.

The one exception is if all non-archimedean primes
have good reduction and the archimedean filled Julia set is an
interval of length $4$.  This occurs for the Chebyshev
polynomial $\phi(z)=z^2 - 2$,
which has filled Julia set $[-2,2]$.
In this special case, after removing the points
$\infty$ and $2$, the rest of the
preperiodic points can be covered by four half-open intervals
of length $1$ at each archimedean prime.
Since there are no non-archimedean bad primes,
we obtain a bound of $2 + 4^D$
for the total number of preperiodic
points in $\PK$.
\end{remark}

\begin{remark}
Another approach to finding a cutoff $\sigma$
which forces Case~1 would be to consider
the set $T\subseteq S$ consisting of non-archimedean
bad primes $v$ at which there are actually $K$-rational
preperiodic points $x,y$ for which $|x-y|_v >1$.
For such primes, the exponent of $-1/[(d-1)(d-2)]$
in Lemma~\ref{lem:minrad} could be improved to $-1/(d-1)$.
Unfortunately, there may not be very many such primes.
As a result, although the exponent of $(d-1)(d-2)$ in
the definition of $\sigma$ could be improved to $(d-1)$,
it would come at the expense of introducing a extra
factor like $\beta^D$ into the formula for $\sigma$.
\end{remark}

\begin{remark}
For large degrees $d$, one can obtain slightly smaller
bounds by using more than one big bad prime $w$.  There is,
of course, a trade-off.  While using $\ell\geq 2$ big primes $w$
ultimately increases the coefficient $A$ of $-N$ in the exponent
of \eqref{eq:midbound}, it also increases the number of
pieces $\{X_k\}$ from $2$ to $2^{\ell}$.
It appears that the
optimal number of such primes to use is $\ell\approx 2\log_2(d-1)$.
The improved bound for the number of rational preperiodic
points would be roughly the old bound divided by $2\log_2(d-1)$,
for large $d$.
However, the proof would be vastly more complicated, especially
dealing with the archimedean primes.  The slight
improvement seems not to be worth the increased difficulty,
especially given that the resulting bound would still be very
far from the conjectured uniform bound.
\end{remark}

We close by presenting a slight strengthening
of Theorem~\ref{thm:global} in the simplest case.

\begin{example}
Let $K=\QQ$ (so that $D=s_{\infty}=1$, and $n_v=1$ for all $v\in M_{\QQ}$)
and $d=2$.  That is, we wish to bound the number of
rational preperiodic points of a quadratic polynomial $\phi\in\QQ[z]$.
It is of
course well known that any such polynomial is conjugate
over $\QQ$ to one of the form $\phi_c(z)=z^2 + c$, with $c\in\QQ$.

Let us suppose that $\phi_c$ has at least one preperiodic point in $\QQ$.
This supposition
implies that $c = j/m^2$ for some relatively prime integers
$j,m\in\ZZ$, and that $-\infty < c\leq 1/4$;
see, for example,
Proposition~6.7 of \cite{CaGol}.
(One can also easily establish that $j$ must satisfy one of approximately
$2^s$ congruences modulo $m$, but we do not need that here.)
For non-archimedean primes $v$ of $\ZZ$, we have $R_v = |m|_v^{-1}$
if $v$ is odd, and $R_2 = \max\{|m/2|_2^{-1}, 1\}$.
(Note that if $4\nmid m$, then $\phi_c$ has good reduction at $v=2$,
after a change of coordinates.)
In addition, for $c<0$, we have $R_{\infty} = (1 + \sqrt{1-4c})/2$.

We will be in Case~1 of the
proof of Theorem~\ref{thm:global} provided there is
some prime $v$ with $R_v\geq 4$.
Still assuming that
there is at least one preperiodic point in $\QQ$,
Lemma~\ref{lem:minrad} says that such a prime must exist
unless the only bad primes are $\infty$, $2$, and $3$,
and $R_2,R_3,R_{\infty}<4$.  By our characterization of $R_v$
above, this means that the denominator $m$ is a divisor of $12$,
and that $-12< c\leq 1/4$.  There are only finitely many rational
numbers of the form $c=j/144$ between $-12$ and $1/4$, and a simple computer
search shows none of the corresponding polynomials $\phi_c$ has
more than eight preperiodic points in $\QQ$.  (For five such
values of $c$, namely
$-21/16$,
$-29/16$,
$-91/36$,
$-133/144$,
and $-1333/144$, there are exactly eight preperiodic points in $\QQ$.
Incidentally, there are infinitely many values $c\in\QQ$
for which $\phi_c$ has at least eight preperiodic point in $\QQ$,
by Theorem~2 of \cite{Poo}.)

For all other $c$, we are in Case~1, which means
$t=s + D\log d/(2\log 2) = s + 1/2$, and $\beta=1$.
If $s=1$, then only the archimedean prime is bad,
and in light of Remark~\ref{rem:improve}, there are
at most five preperiodic points in $\QQ$; in fact,
there are at most four for $s=1$ and $c\neq -2$.
The only remaining possibility is that $s\geq 2$ and $\beta=1$,
in which case the
number of preperiodic points in $\QQ$ is at most
$$(2s + 1)\left[ \log_2(2s + 1) + \log_2( \log_2(2s+1) - 1) + 2\right].$$
Since this bound is greater than eight
even for $s=2$, it holds even without making the exceptions
from the previous paragraph.
\end{example}

\bibliographystyle{plain}

\end{document}